\documentclass{amsproc}

\usepackage[arrow, curve, frame, matrix, graph, tips]{xy}
\usepackage{amssymb}
\usepackage{amsthm}
\usepackage{mtens}

\newtheorem{theorem}{Theorem}[section]
\newtheorem{lemma}[theorem]{Lemma}
\newtheorem{proposition}[theorem]{Proposition}
\newtheorem{corollary}[theorem]{Corollary}

\theoremstyle{definition}
\newtheorem{definition}[theorem]{Definition}
\newtheorem{example}[theorem]{Example}
\newtheorem{examples}[theorem]{Examples}
\newtheorem{non-example}[theorem]{Non-Example}

\theoremstyle{remark}
\newtheorem{remark}[theorem]{Remark}

\newcommand{\TriTwoCell}[7]{\xymatrix{
{#1} \ar[rr]^-{#5} \ar[dr]_{#4} \save \POS?="dom" \restore
&& {#2} \ar[dl]^{#6} \save \POS?="cod" \restore \\ & {#3}
\POS "dom"; "cod" **@{} ?(.35) \ar@{=>}^{#7} ?(.65)
}}

\newcommand{\TwoDiagRel}[3]{\begin{xy}
(0,0);<2em,0em>:<0em,2em>::
(-1,0)*-!R{\xybox{#1}};
(1,0)*-!L{\xybox{#3}}
**@{} ?*{#2}
\end{xy}}

\newcommand{\LaxSq}[9]{\xymatrix{{#1} \ar[r]^-{#6}
\ar[d]_{#5} \save \POS?="dom" \restore
& {#2} \ar[d]^{#7} \save \POS?="cod" \restore \\
{#3} \ar[r]_-{#8} & {#4}
\POS "dom"; "cod" **@{} ?(.35) \ar@{=>}^{#9} ?(.65)}}

\newcommand{\MNDmorphism}[6]{\xymatrix{{#1} \ar[r]^-{#5} \ar[d]_{#2}
& {#3} \save \POS="dom" \restore \ar[d]^{#4} \\
{#1} \save \POS="cod" \restore \ar[r]_-{#5} & {#3}
\POS "dom"; "cod" **@{} ?(.35) \ar@{=>}^{#6} ?(.65)}}

\newcommand{\PbSq}[8]{\xymatrix{{#1} \ar[d]_{#5}
\save \POS?(.3)="lpb" \restore
\ar[r]^-{#6} \save \POS?(.3)="tpb" \restore
& {#2} \ar[d]^{#7} \save \POS?(.3)="rpb" \restore \\
{#4} \ar[r]_-{#8} \save \POS?(.3)="bpb" \restore & {#3}
\POS "rpb"; "lpb" **@{}; ?!{"bpb";"tpb"}="cpb" **@{}; ? **@{-};
"tpb"; "cpb" **@{}; ? **@{-}}}

\newcommand{\rel}{\SelectTips{cm}{}\object@{/}}

\begin{document}

\title{Algebras of higher operads as enriched categories}

\author{Michael Batanin}
\address{Department of Mathematics,
Macquarie University}
%\curraddr{}
\email{mbatanin@ics.mq.edu.au}
\thanks{}
% author two information
\author{Mark Weber}
\address{Laboratoire PPS, Universit\'e Paris Diderot -- Paris 7}
%\curraddr{}
\email{weber@pps.jussieu.fr}
\thanks{}
\maketitle
%%%%%%%%%%%%%%%%%%%%%%%%%%%%%%%%%%%%%%%%%%%%%%%%%%
%%%%%%%%%%%%%%%%%%%%%%%%%%%%%%%%%%%%%%%%%%%%%%%%%%
\begin{abstract}
We decribe the correspondence between normalised $\omega$-operads in the sense of \cite{Bat98} and certain lax monoidal structures on the category of globular sets. As with ordinary monoidal categories, one has a notion of category enriched in a lax monoidal category. Within the aforementioned correspondence, we provide also an equivalence between the algebras of a given normalised $\omega$-operad, and categories enriched in globular sets for the induced lax monoidal structure. This is an important step in reconciling the globular and simplicial approaches to higher category theory, because in the simplicial approaches one proceeds inductively following the idea that a weak $(n+1)$-category is something like a category enriched in weak $n$-categories, and in this paper we begin to reveal how such an intuition may be formulated in terms of the machinery of globular operads.
\end{abstract}
%%%%%%%%%%%%%%%%%%%%%%%%%%%%%%%%%%%%%%%%%%%%%%%%%%
%%%%%%%%%%%%%%%%%%%%%%%%%%%%%%%%%%%%%%%%%%%%%%%%%%
\section{Introduction}

The subject of enriched category theory \cite{Kel82} was brought to maturity by the efforts of Max Kelly and his collaborators. Max also had a hand in the genesis of the study of operads, and in \cite{KelOp} which for a long time went unpublished, he layed the categorical basis for their further analysis. It is with great pleasure that we are able to present the following paper, which relates enriched category theory and the study of higher operads, in dedication to a great mathematician.

In the combinatorial approach to defining and working with higher categorical structures, one uses globular operads to say what the structures of interest are in one go. However in the simplicial approaches to higher category theory, one proceeds inductively following the idea that a weak $(n+1)$-category is something like a category enriched in weak $n$-categories.
This is the first in a series of papers whose purpose is to reveal and study the inductive aspects hidden within the globular operadic approach.

An $\omega$-operad in the sense of \cite{Bat98} can be succinctly described as a cartesian monad morphism $\alpha:A{\rightarrow}\ca T$, where $\ca T$ is the monad on the category $\PSh {\G}$ of globular sets whose algebras are strict $\omega$-categories. The algebras of the given operad are just the algebras of the monad $A$. Among the $\omega$-operads, one can distinguish the \emph{normalised} ones, which don't provide any structure at the object level, so that one may regard a globular set $X$ and the globular set $AX$ as having the same objects. For example, the operad constructed in \cite{Bat98} to define weak-$\omega$-categories, and indeed any $\omega$-operad that has been constructed to give a definition of weak-$\omega$-category, is normalised. One of the main results of this paper, corollary(\ref{cor:3way-equiv}), provides two alternative views of normalised operads: as $M\ca T$-operads and as $\ca T$-multitensors.

The notion of $T$-operad, and more generally of $T$-multicategory, makes sense for any cartesian monad $T$ on a finitely complete category $\ca V$ (see \cite{Lei}). A $T$-operad can be defined as a cartesian monad morphism into $T$, in the same way as we have already outlined in the case $T=\ca T$ above. Under certain conditions on $\ca V$ and $T$, one has a monad $M$ on $\ca V$ which is also cartesian and whose algebras are monoids in $\ca V$, this monad distributes with $T$, and the composite monad $MT$ is also cartesian, so one can consider $MT$-operads. All of this is so in the case $T=\ca T$.

On the other hand a \emph{multitensor} structure on a category $\ca V$ is just another name for the structure of a lax monoidal category on $\ca V$. This general notion has been discussed both in \cite{DS03} within the framework of lax monoids, and in \cite{Bat02} where it is expressed in the language of internal operads. A multitensor is like a monoidal structure, except that the coherences are not necessarily invertible, and one works in an ``unbiased'' setting defining an $n$-ary tensor product for all $n \in \N$. Just as with monoidal categories one can consider categories enriched in a lax monoidal category. In particular if $\ca V$ has cartesian products and $T$ is a monad on $\ca V$, one can define a canonical multitensor $T^{\times}$ on $\ca V$, with the property that categories enriched in $(\ca V,T^{\times})$ are exactly categories enriched in $\Alg T$ for the cartesian tensor product. When $\ca V$ is lextensive and $T$ is a p.r.a monad in the sense of \cite{Fam2fun}, one can define a $T$-multitensor in an analogous way to the definition of $T$-operad: as a cartesian multitensor morphism into $T^{\times}$. These assumptions on $\ca V$ and $T$ are a little stronger than asking that $T$ be a cartesian monad, and are clearly satisfied for all examples of interest for us such as $T=\ca T$.

The correspondence between normalised $\ca T$-operads and $\ca T$-multitensors already discussed also includes an important feature at the level of algebras. Namely the algebras of a given normalised $\ca T$-operad $\alpha:A{\rightarrow}\ca T$ correspond to categories enriched in the associated $\ca T$-multitensor. In this way, any higher categorical structure definable by a normalised $\ca T$-operad is expressed as a category enriched in $\PSh {\G}$ for a canonically defined lax monoidal structure on $\PSh {\G}$.

This paper is organised as follows. In section(\ref{sec:lax-mon}) we recall the definition of a lax monoidal category and of categories enriched therein, and give the example of $T^{\times}$. Multitensors, that is lax monoidal structures, generalise non-symmetric operads, and sections(\ref{ssec:mult-mon}) and (\ref{ssec:mult-fmon}) explain how basic operad theory generalises to multitensors. In section(\ref{ssec:mult-mon}) we see how under certain conditions, one may regard multitensors as monoids for a certain monoidal structure, which generalises the substitution tensor product of collections familiar from the theory of operads. Proposition(\ref{prop:DistV}) is in fact a special case of proposition(2.1) of \cite{DS03}. Nevertheless we give a self-contained account of proposition(\ref{prop:DistV}) and related notions, to keep the exposition relatively self-contained and as elementary as possible for our purposes. In section(\ref{ssec:mult-fmon}) we explain how one can induc!
 
 e a monad from a multitensor. The theory of $T$-multitensors, which is the multitensorial analogue of the theory of $T$-operads described in \cite{Lei}, is given in section(\ref{ssec:masop}), and it is at this level of generality that one sees the equivalence between $\ca T$-multitensors and $M\ca T$-operads.

From this point in the paper we begin working directly with the case $T=\ca T$. In section(\ref{sec:wcat-monad}) we give a self-contained inductive description of the monad $\ca T$. This is a very beautiful mathematical object. It is a p.r.a monad and its functor part preserves coproducts. It has another crucial property, called \emph{tightness}, which implies that for any endofunctor $A$, if a cartesian transformation $\alpha:A{\rightarrow}\ca T$ exists then it is unique. This property is very useful, for instance when building up a description of $\ca T$ one need not check the monad axioms because these come for free once one has given cartesian transformations $\eta:1{\rightarrow}\ca T$ and $\mu:\ca T^2{\rightarrow}\ca T$. The inductive description of $\ca T$ given here is closely related to the wreath product of Clemens Berger \cite{Ber07}.

In section(\ref{ssec:corr}) we give the correspondence between normalised $\ca T$-operads and $\ca T$-multitensors, as well as the identification between the algebras of a given normalised operad and categories enriched in the associated multitensor. In the final section we explain how our results may be adapted to normalised $n$-operads, that is to finite dimensions, and then explain how the algebras of $\ca T$, which we defined as a combinatorial object, really are strict $\omega$-categories. This last fact is of course well-known, but the simplicity and canonicity of our proof is a pleasant illustration of the theory developed in this paper.

The work discussed here is in a sense purely formal. Everything works at a high level of generality. Things become more interesting and subtle when we wish to lift the lax monoidal structures we obtain on $\PSh {\G}$, or one of its finite dimensional analogues, to the category of algebras of another operad. For example already in this paper, one can see that the lax monoidal structure $\ca T^{\times}$ on $\PSh {\G}$ corresponds to cartesian product of $\ca T$-algebras, in the sense that they give the same enriched categories. It is from the general theory of such lifted lax monoidal structures that the Gray tensor product and its many variants, and many other examples, will be captured within our framework. These issues will be the subject of \cite{EnHopII}.

%%%%%%%%%%%%%%%%%%%%%%%%%%%%%%%%%%%%%%%%%%%%%%%%%%
%%%%%%%%%%%%%%%%%%%%%%%%%%%%%%%%%%%%%%%%%%%%%%%%%%
\section{Lax Monoidal Categories}\label{sec:lax-mon}

In this section we recall the notion of \emph{lax monoidal category}, which is a generalisation of the well-known concept of monoidal category. As with monoidal categories, one can consider categories enriched in a lax monoidal category. Any monad $T$ on a finitely complete category $\ca V$ defines a canonical lax monoidal structure $T^{\times}$ on $\ca V$, and for this structure enriched categories correspond to categories enriched in $\Alg T$ regarded as monoidal via cartesian product.

Given a 2-monad $T$ on a 2-category $\ca K$ one may consider
\emph{lax algebras} for $T$. A lax $T$-algebra structure on an object $A \in \ca K$ is a triple $(a,u,\sigma)$ consisting of an action $a:TA{\rightarrow}A$ together with 2-cells
%%%%%
\[ \TwoDiagRel
{\xymatrix{A \ar[rr]^-{{\eta}_{A}} \ar[dr]_{1_A} \save \POS?="codlam" \restore
&& {{T}A} \ar[dl]^{a} \save \POS?="domboth" \restore \\ & A
\POS "codlam";"domboth"**@{} ?  /-.5em/ \ar@{=>}^-{u} ?  /.5em/ \restore}}
{}
{\xymatrix{{T^2A} \ar[r]^-{{\mu}_{A}} \ar[d]_{{T}a} \save \POS?="domass" \restore
& {{T}A} \ar[d]^{a} \save \POS?="codass" \restore \\ {{T}A} \ar[r]_-{a} & A
\save \POS"domass";"codass"**@{} ?  /-.5em/ \ar@{=>}^-{\sigma} ?  /.5em/ \restore}} \]
%%%%%
satisfying some well-known axioms. See \cite{LkCodesc} for a complete description
of these axioms, and of the 2-category $\LaxAlg T$. When $T$ is the identity, lax algebras are just monads in $\ca K$. The example most important for us however is when $T$ is the monoid monad $\ca M$ on $\CAT$.
\begin{definition}\label{def:multitensor}
A \emph{multitensor} on a category $\ca V$ is a lax $\ca M$-algebra structure $(E,u,\sigma)$ on $\ca V$.
A category $\ca V$ equipped with a multitensor structure is called a \emph{lax monoidal category}. When $u$ is the identity the multitensor and lax monoidal structure are said to be \emph{normal}.
\end{definition}
\noindent We shall now unpack this definition. Since $\ca M\ca V = \coprod\limits_{n{\geq}0} \ca V^n$
a functor \[ E:\ca M\ca V{\rightarrow}\ca V \] amounts to functors $E_n:\ca V^n{\rightarrow}\ca V$ for $n \in \N$.

Before proceeding further we digress a little on notation. For functors of many variables we shall use some space saving notation: we deem that the following expressions
\[ \begin{array}{lcccr} E_n(X_1,...,X_n) && \opE\limits_{1{\leq}i{\leq}n} X_i &&  \opE\limits_iX_i
\end{array} \]
are synonymous, and we will frequently use the latter, often leaving the ``$n$'' unmentioned when no confusion would result. In particular for $X{\in}\ca V$ and $1{\leq}i{\leq}n$, $\opE\limits_iX$ denotes $E_n(X,...,X)$. We identify the number $n$ with the ordered set $\{1,...,n\}$ and we refer to elements of the ordinal sum $n_{\bullet}:=n_1{+}...{+}n_k$ as pairs $(i,j)$ where $1{\leq}i{\leq}k$ and $1{\leq}j{\leq}n_i$. Following these conventions $\opE\limits_i\opE\limits_jX_{ij}$ and $\opE\limits_{ij}X_{ij}$ are synonymous with
\[ E_k(E_{n_1}(X_{11},...,X_{1n_1}),...,E_{n_k}(X_{k1},...,X_{kn_k})) \]
and
\[ E_{n_{\bullet}}(X_{11},...,X_{1n_1},......,X_{k1},...,X_{kn_k}) \]
respectively. We will use multiply indexed expressions (like $\opE\limits_i\opE\limits_{jk}\opE\limits_lX_{ijkl}$) to more efficiently convey expressions that have multiple layers of brackets and applications of $E$'s.

The remaining data for a multitensor on $\ca V$ amounts to maps
\[ \begin{array}{lcr}
u_{X}:X{\rightarrow}E_1X &&
\sigma_{X_{ij}}:\opE\limits_i\opE\limits_jX_{ij}{\rightarrow}\opE\limits_{ij}X_{ij} \end{array} \]
%%%%%
that are natural in the arguments and satisfy
%%%%%%%%
\[ \xy (0,0); (10,0):
(0,0)*{\xybox{\xymatrix @C=1em {{\opE\limits_iX_i} \ar[r]^-{u\opE\limits_i} \ar[d]_{1} \save \POS?(.4)="domeq" \restore & {E_1\opE\limits_iX_i} \ar[dl]^{\sigma} \save \POS?(.4)="codeq" \restore \\ {\opE\limits_iX_i} \POS "domeq"; "codeq" **@{}; ?*{=}}}};
(4,0)*{\xybox{\xymatrix @C=1.5em {{\opE\limits_i\opE\limits_j\opE\limits_kX_{ijk}} \ar[r]^-{{\sigma}\opE\limits_k} \ar[d]_{\opE\limits_i\sigma} \save \POS?="domeq" \restore & {\opE\limits_{ij}\opE\limits_kX_{ijk}} \ar[d]^{\sigma} \save \POS?="codeq" \restore \\
{\opE\limits_i\opE\limits_{jk}X_{ijk}} \ar[r]_-{\sigma} & {\opE\limits_{ijk}X_{ijk}} \POS "domeq"; "codeq" **@{}; ?*{=}}}};
(8,0)*{\xybox{{\xymatrix @C=1em {{\opE\limits_iE_1X_i} \ar[dr]_{\sigma} \save \POS?(.4)="domeq" \restore & {\opE\limits_iX_i} \ar[d]^{1} \save \POS?(.4)="codeq" \restore  \ar[l]_-{\opE\limits_iu} \\ & {\opE\limits_iX_i}} \POS "domeq"; "codeq" **@{}; ?*{=}}}}
\endxy \]
%%%%%%%%
Thus a multitensor is very much like a functor-operad in the sense of \cite{McClureSmith}, except that there are no symmetric group actions with respect to which the substitutions are equivariant{\footnotemark{\footnotetext{
More precisely, functor-operads in the sense of \cite{McClureSmith} are normal lax algebras for the \emph{symmetric} monoidal category 2-monad on $\CAT$.}}}. An equivalent formulation of definition(\ref{def:multitensor}), in the language of \cite{Bat02}, is that a multitensor on $\ca V$ is a non-symmetric operad internal to the endomorphism operad of $\ca V$. 
\begin{example}\label{ex:mult-prod}
A normal multitensor on $\ca V$ such that $\sigma$ is invertible is just a monoidal structure on $\ca V$, with $E_n$ playing the role of the $n$-fold tensor product. In the case where $\ca V$ is finitely complete and $E_n$ is $n$-fold cartesian product and for the sake of the next example, we denote the isomorphism ``$\sigma$'' as
\[ \begin{array}{c} {\iota : \prod\limits_i\prod\limits_j X_{ij} \rightarrow \prod\limits_{ij} X_{ij}} \end{array} \]
\end{example}
\begin{example}\label{ex:mult-Tcross}
Let $T$ be a monad on a finitely complete category $\ca V$. Denote by
\[ \begin{array}{c} {k_{X_i} : T \prod\limits_iX_i \rightarrow \prod\limits_iTX_i} \end{array} \]
the canonical maps which measure the extent to which $T$ preserves products. One defines a multitensor $(T^{\times},u,\sigma)$ as follows:
\[ \begin{array}{c} {T^{\times}_k(X_1,...,X_k) = \prod\limits_{1{\leq}i{\leq}n} T(X_i)} \end{array} \]
$u$ is the unit $\eta_X:X{\rightarrow}TX$ of the monad, and $\sigma$ is defined as the composite
\[ \xymatrix{{\prod\limits_iT\prod\limits_jTX_{ij}} \ar[r]^-{\prod\limits_ikT} & {\prod\limits_i\prod\limits_jT^2X_{ij}} \ar[r]^-{\iota\mu} & {\prod\limits_{ij}TX_{ij}}} \]
\end{example}
\noindent For the remainder of this section let $(\ca V,E)$ be a lax monoidal category.
\begin{definition}
An \emph{$E$-category} $(X,\kappa)$, or in other words a category enriched in $(\ca V,E)$, consists of
\begin{itemize}
\item  a set $X_0$ of objects.
\item  for all pairs $(x_0,x_1)$ of elements of $X_0$, an object $X(x_0,x_1)$ of $\ca V$. These objects are called the \emph{homs} of $X$.
\item  for all $n \in \N$ and $(n{+}1)$-tuples $(x_0,...,x_n)$ of elements of $X_0$, maps
\[ \kappa_{x_i} : \opE\limits_{1{\leq}i{\leq}n} X(x_{i-1},x_{i}) \rightarrow X(x_0,x_n) \]
called the \emph{compositions} of $X$.
\end{itemize}
satisfying unit and associative laws, which say that
%%%%%%%%
\[ \xygraph{
{\xybox{\xygraph{!{0;(2,0):(0,.6)::}
{X(x_0,x_1)} (:[r]{E_1X(x_0,x_1)}^-{u} :[d]{X(x_0,x_1)}="bot"^{\kappa}, :"bot"_{\id})}}}
[r(5)][d(.15)]
{\xybox{\xygraph{!{0;(2.75,0):(0,.5)::}
{\opE\limits_i\opE\limits_jX(x_{(ij)-1},x_{ij})} (:[r]{\opE\limits_{ij}X(x_{(ij)-1},x_{ij})}^-{\sigma} :[d]{X(x_0,x_{mn_m})}="bot"^{\kappa},:[d]{\opE\limits_iX(x_{(i1)-1},x_{in_i})}_{\opE\limits_i\kappa} :"bot"_-{\kappa})}}}} \]
%%%%%%%%
commute, where $1{\leq}i{\leq}m$, $1{\leq}j{\leq}n_i$ and $x_{(11)-1}{=}x_0$. Since a choice of $i$ and $j$ references an element of the ordinal $n_{\bullet}$, the predecessor $(ij){-}1$ of the pair $(ij)$ is well-defined when $i$ and $j$ are not both $1$. An \emph{$E$-monoid} is an $E$-category with one object.
\end{definition}
\begin{definition}
Let $(X,\kappa)$ and $(Y,\lambda)$ be $E$-categories. An \emph{$E$-functor} $f:(X,\kappa){\rightarrow}(Y,\lambda)$ consists of a function $f_0:X_0{\rightarrow}Y_0$, and for all pairs $(x_0,x_1)$ from $X_0$, arrows
\[ f_{x_0,x_1} : X(x_0,x_1) \rightarrow Y(fx_0,fx_1) \]
satisfying a functoriality axiom, which says that
\[ \xymatrix{{\opE\limits_iX(x_{i-1},x_i)} \ar[r]^-{\opE\limits_if} \ar[d]_{\kappa} & {\opE\limits_iY(fx_{i-1},fx_i)} \ar[d]^{\lambda} \\ {X(x_0,x_n)} \ar[r]_-{f} & {Y(fx_0,fx_n)}} \]
commutes. We denote by $\Enrich E$ the category of $E$-categories and $E$-functors, and by $\Mon E$ the full subcategory of $\Enrich E$ consisting of the $E$-monoids.
\end{definition}
\begin{example}\label{ex:nso-as-mt}
A non-symmetric operad
\[ \begin{array}{lcr} {(A_n \, \, : \, \, n \in \N)} & {u:I \rightarrow A_1} &
{\sigma:A_k \tensor A_{n_1} \tensor ... \tensor A_{n_k} \rightarrow A_{n_{\bullet}}} \end{array} \]
in a braided monoidal category $\ca V$ defines a multitensor $E$ on $\ca V$ via the formula
\[ \opE\limits_{1{\leq}i{\leq}n} X_i = A_n \tensor X_1 \tensor ... \tensor X_n \]
with $u$ and $\sigma$ providing the structure maps in the obvious way. The category $\Mon E$ of $E$-monoids is the usual category of algebras of $A$, and thus $E$-categories are a natural notion of ``many object algebra'' for an operad $A$.
\end{example}
\noindent Our notation for multitensors makes evident the analogy with monads and algebras: a multitensor $E$ is analogous to a monad and an $E$-category is the analogue of an algebra for $E$. In particular observe that the following basic facts are instances of the axioms for the lax monoidal category $(\ca V,E)$ and categories enriched therein.
\begin{lemma}\label{lem:basic-obs}
\begin{enumerate}
\item  $(E_1,u,\sigma)$ is a monad on $\ca V$.
\item  The monad $E_1$ acts on $E_n$ for all $n \in \N$, that is
\[ \sigma : E_1\opE\limits_iX_i \rightarrow \opE\limits_iX_i \]
is an $E_1$-algebra structure on $\opE\limits_iX_i$.\label{E1-action}
\item  With respect to the $E_1$-algebra structures of (\ref{E1-action}) all of the components of $\sigma$ are $E_1$-algebra morphisms.
\item  Each hom of an $E$-category $(X,\kappa)$ is an $E_1$-algebra, with the algebra structure on $X(x_0,x_1)$ given by \[ \kappa : E_1X(x_0,x_1) \rightarrow X(x_0,x_1). \]
\label{E1-action2}
\item  With respect to the $E_1$-algebras of (\ref{E1-action}) and (\ref{E1-action2}), all the components of $\kappa$ are morphisms of $E_1$-algebras.
\end{enumerate}
\end{lemma}
\begin{proposition}\label{prop:Tcross-cats}
Let $T$ be a monad on a finitely complete category $\ca V$. Regarding $\Alg T$ as a monoidal category via cartesian product one has \[ \Enrich {T^{\times}} \iso \Enrich {(\Alg T)} \] commuting with the forgetful functors into $\Set$.
\end{proposition}
\begin{proof}
Let $X_0$ be a set and for $a,b \in X_0$ let $X(a,b) \in \ca V$. Suppose that
\[ \begin{array}{c} {\kappa_{x_i}:\prod\limits_{i}TX(x_{i-1},x_i){\rightarrow}X(x_0,x_n)} \end{array} \]
for each $n \in \N$ and $x_0,...,x_n$ in $X_0$, are the structure maps for a $T^{\times}$-category structure. Then by lemma(\ref{lem:basic-obs}) the $\kappa_{a,b}:TX(a,b){\rightarrow}X(a,b)$ are algebra structures for the homs, and for $x_{ij} \in X_0$ with $1{\leq}i{\leq}k$ and $1{\leq}j{\leq}n_i$ one has the inner regions of
%%%%%%%
\[ \xygraph{!{0;(3.7,0):(0,.4)::}
{T\prod\limits_iX(x_{i-1},x_i)}="t1" ([d] {\prod\limits_iTX(x_{i-1},x_i)}="m1" [d] {\prod\limits_iX(x_{i-1},x_i)}="b1",[r] {T\prod\limits_iTX(x_{i-1},x_i)}="t2" ([r] {TX(x_0,x_n)}="t3" [dd] {X(x_0,x_n)}="b3" [l(.5)u(.6)] {\prod\limits_iTX(x_{i-1},x_i)}="m3",[l(.25)d] {\prod\limits_iT^2X(x_{i-1},x_i)}="m2" [r(.15)d] {\prod\limits_iTX(x_{i-1},x_i)}="b2"))
"t1":"t2"^-{T\prod\limits_i\eta} :"t3"^-{T\kappa} :"b3"^-{\kappa} "t1":"m1"_-{k} :"b1"_-{\prod\limits_i\kappa} :"b2"_-{\prod\limits_i\eta} :"b3"_-{\kappa} "t2":"m2"^-{k} :"b2"^-{\prod\limits_iT\kappa} "m1":"m2"^-{\prod\limits_i\eta{T}} :"m3"^-{\prod\limits_i\mu} :"b3"^-{\kappa}} \]
%%%%%%%
commutative, and the commutative outer region is the associativity axiom for the composites
\[ \xymatrix{{\kappa'_{x_i}:\prod\limits_iX(x_{i-1},x_i)} \ar[r]^-{\prod\limits_i\eta} & {\prod_iTX(x_{i-1},x_i)} \ar[r]^-{\kappa_{x_i}} & {X(x_0,x_n)}} \]
for each $x_0,...,x_n \in X_0$. Taking the product structure on $\Alg T$ as normal, the unit axiom for the $\kappa'$ is clearly satisfied, and so they are the structure maps for a $(\Alg T)$-category structure. Conversely given algebra structures $\kappa_{a,b}$ and structure maps $\kappa'_{x_i}$ one can define $\kappa_{x_i}$ as the composite
\[ \xymatrix{{\prod\limits_iTX(x_{i-1},x_i)} \ar[rr]^-{\prod\limits_i\kappa_{x_{i-1},x_i}} && {\prod\limits_iX(x_{i-1},x_i)} \ar[r]^-{\kappa'_{x_i}} & {X(x_0,x_n)}} \]
and since the regions of
%%%%%%%
\[ \xygraph{!{0;(3.25,0):(0,.4)::}
{\prod\limits_iT\prod\limits_jTX(x_{(ij)-1},x_{ij})}="t1" [r] {\prod\limits_{ij}T^2X(x_{(ij)-1},x_{ij})}="t2" [r] {\prod\limits_{ij}TX(x_{(ij)-1},x_{ij})}="t3" [d] {\prod\limits_{ij}X(x_{(ij)-1},x_{ij})}="m3" [d] {X(x_0,x_n)}="b3" [l] {\prod\limits_iX(x_{i-1},x_i)}="b2" [l] {\prod\limits_iTX(x_{i-1},x_i)}="b1" [u] {\prod\limits_iT\prod\limits_jX(x_{(ij)-1},x_{ij})}="m1" [r] {\prod\limits_{ij}T^2(x_{(ij)-1},x_{ij})}="m2"
"t1":"t2"^-{k} :"t3"^-{\prod\limits_{ij}\mu} :"m3"^-{\prod\limits_{ij}\kappa} :"b3"^-{\kappa'} "t1":"m1"_-{\prod\limits_iT\prod\limits_j\kappa} :"b1"_-{\prod\limits_iT\kappa'} :"b2"_-{\prod\limits_i\kappa} :"b3"_-{\kappa'} "m1":"m2"^-{k} :"m3"^-{\prod\limits_{ij}\kappa} :"b2"_{\prod\limits_i\kappa'} "t2":"m2"^{\prod\limits_{ij}T\kappa}} \]
%%%%%%%
commute, the commutativity of the outside of this diagram shows that the $\kappa_{x_i}$ satisfy the associativity condition of a $T^{\times}$-category structure, and the unit axiom follows from the unit $T$-algebra axiom on the homs. The correspondence just described is clearly a bijection, and completes the description of the isomorphism on objects over $\Set$. 

Let $f_0:X_0{\rightarrow}Y_0$ be a function,
\[ \begin{array}{lccr} {\kappa_{x_i}:\prod\limits_{i}TX(x_{i-1},x_i){\rightarrow}X(x_0,x_n)} &&
{\lambda_{y_i}:\prod\limits_{i}TY(y_{i-1},y_i){\rightarrow}Y(y_0,y_n)} \end{array} \]
be the structure maps for $T^{\times}$-categories $X$ and $Y$, $\kappa'_{x_i}$ and $\lambda'_{y_i}$ be the associated $(\Alg T)$-category structures, and
\[ f_{a,b} : X(a,b) \rightarrow Y(fa,fb) \]
for $a,b \in X_0$ be maps in $\ca V$. In the following display the diagram on the left
\[ \TwoDiagRel {\xymatrix{{\prod\limits_iX(x_{i-1},x_i)} \ar[r]^-{\prod\limits_if} \ar[d]_{\prod\limits_i\eta} & {\prod\limits_iY(y_{i-1},y_i)} \ar[d]^{\prod\limits_i\eta} \\
{\prod\limits_iTX(x_{i-1},x_i)} \ar[r]^-{\prod\limits_iTf} \ar[d]_{\kappa} & {\prod\limits_iTY(y_{i-1},y_i)} \ar[d]^{\lambda} \\ {X(x_0,x_n)} \ar[r]_-{f_{x_0,x_n}} & {Y(y_0,y_n)}}} {}
{\xymatrix{{\prod\limits_iTX(x_{i-1},x_i)} \ar[r]^-{\prod\limits_iTf} \ar[d]_{\prod\limits_i\kappa} & {\prod\limits_iTY(y_{i-1},y_i)} \ar[d]^{\prod\limits_i\lambda} \\
{\prod\limits_iX(x_{i-1},x_i)} \ar[r]^-{\prod\limits_if} \ar[d]_{\kappa'} & {\prod\limits_iY(y_{i-1},y_i)} \ar[d]^{\lambda'} \\ {X(x_0,x_n)} \ar[r]_-{f_{x_0,x_n}} & {Y(y_0,y_n)}}} \]
explains how the $T^{\times}$-functor axiom for the $f_{a,b}$ implies the $(\Alg T)$-functor axiom, and the diagram on the right shows the converse.
\end{proof}

%%%%%%%%%%%%%%%%%%%%%%%%%%%%%%%%%%%%%%%%%%%%%%%%%%
%%%%%%%%%%%%%%%%%%%%%%%%%%%%%%%%%%%%%%%%%%%%%%%%%%
\section{Distributive multitensors as monoids}\label{ssec:mult-mon}

It is well-known that monads on a category $\ca V$ are monoids in the strict monoidal category $\End(\ca V)$ of endofunctors of $\ca V$ whose tensor product is given by composition. Given the analogy between monads and multitensors, one is led to ask under what circumstances are multitensors monoids in a certain monoidal category. One natural answer to this question, that we shall present now, requires that we restrict attention to \emph{distributive} multitensors to be defined below. Throughout this section $\ca V$ is assumed to have coproducts.
\begin{definition}\label{def:distMF}
A functor $E:\ca M\ca V{\rightarrow}\ca V$ is \emph{distributive} when for all $n \in \N$, $E_n$ preserves coproducts in each variable. We denote by $\Dist(\ca V)$ the category whose objects are such functors $\ca M\ca V{\rightarrow}\ca V$, and whose morphisms are natural transformations between them. A multitensor $(E,u,\sigma)$ (resp. lax monoidal category $(\ca V,E)$) is said to be \emph{distributive} when $E$ is distributive.
\end{definition}
\begin{examples}\label{ex:Tcross-dist}
In the case where $(\ca V,\tensor,I)$ is a genuine monoidal category, $\ca V$ is distributive in the above sense iff $(X{\tensor}-)$ and $(-{\tensor}X)$ preserve coproducts for each $X \in \ca V$. If in addition $\tensor$ is just cartesian product and $T$ is a monad on $\ca V$ whose functor part preserves coproducts, then the multitensor $T^{\times}$ of example(\ref{ex:mult-Tcross}) is also distributive.
\end{examples}
\noindent When $E$ is distributive we have
\[ \opE_{1{\leq}i{\leq}n} \coprod_{j{\in}J_i} X_{ij} \iso
\coprod_{j_1{\in}J_1} ... \coprod_{j_n{\in}J_n} \opE_iX_{ij_i} \]
for any doubly indexed family $X_{ij}$ of objects of $\ca V$. To characterise distributivity via this formula we must be more precise and say that a certain canonical map between these objects is an isomorphism. It is however more convenient to express all this in terms of coproduct cocones. To state such an equation we must have for each $1{\leq}i{\leq}n$ a family of maps
\[ (c_{ij} : X_{ij} \rightarrow X_{i\bullet} : j \in J_i) \]
which forms a coproduct cocone in $\ca V$. Given a choice for each $i$ of $j \in J_i$, one obtains a map
\[ \opE\limits_ic_{ij} : \opE\limits_i X_{ij} \rightarrow \opE\limits_iX_{i\bullet}, \]
and distributivity says that all such maps together form a coproduct cocone. The morphisms that comprise this cocone are indexed by elements of $\prod\limits_iJ_i$ in agreement with the right hand side of the above formula. For what will soon follow it is worth recalling that the (obviously true) statement ``a coproduct of coproducts is a coproduct'' can be described in a similar way. That is, given $c_{ij}$ as above together with another coproduct cocone
\[ (c_i : X_{i\bullet} \rightarrow X_{\bullet\bullet} : 1{\leq}i{\leq}n), \]
for each choice of $i$ and $j$ one obtains a composite arrow
\[ \xymatrix{{X_{ij}} \ar[r]^-{c_{ij}} & {X_{i\bullet}} \ar[r]^-{c_i} & {X_{\bullet\bullet}}}, \]
and the collection of all such composites is a coproduct cocone.

Define the unit $I$ of $\Dist(\ca V)$ by $I_1{=}1_{\ca V}$ and for $n{\neq}1$, $I_n$ is constant at $\emptyset$. The tensor product $E{\comp}F$ of $E$ and $F$ in $\Dist(\ca V)$ is defined as:
%%%%%%%%
\[  (E \comp F)_n = \coprod_{k{\geq}0}\coprod_{n_1+...+n_k{=}n} \opE\limits_iF_{n_i} \]
%%%%%%%%
and so for all $k$ and $n_i \in \N$ where $1{\leq}i{\leq}k$ we have maps
\[ \xymatrix{{\opE\limits_{i}\opF\limits_{j}} \ar[r]^-{\opc\limits_{ij}} & {\opEoF\limits_{ij}}} \]
which we shall also denote by $c_{(n_1,...,n_k)}$ as convenience dictates. For all $n \in \N$ the set of all such maps such that $n_{\bullet}{=}n$ form a coproduct cocone. In the case where $E{=}I$ one has $\opI\limits_{i}\opF\limits_{j}{\iso}\emptyset$ when $k{\neq}1$, and so
\[ c_{(n)} : F_n \rightarrow (I \comp F)_n \]
is invertible, the inverse of which we denote by $\lambda$. In the case where $F{=}I$ one has $\opE\limits_{i}\opI\limits_{j}{\iso}\emptyset$ when not all the $n_i$'s are $1$, and so
\[ c_{(1,...,1)} : E_n \rightarrow (E \comp 1)_n \]
is invertible, the inverse of which we denote by $\rho$. Given $E$, $F$ and $G$ in $\Dist(\ca V)$, one has for all $r \in \N$, $m_i \in \N$ such that $1{\leq}i{\leq}r$, and $n_{ij} \in \N$ for all $i$ and $1{\leq}j{\leq}m_i$, a composite
\[ \xymatrix{{\opE\limits_i\opF\limits_j\opG\limits_k} \ar[r]^-{\opc\limits_{ij}\opG\limits_k} & {(\opEoF\limits_{ij})\opG\limits_k} \ar[r]^-{\opc\limits_{(ij)k}} & {((E{\comp}F){\comp}G)_{n_{\bullet\bullet}}}} \]
and for all $n \in \N$, the set of all such composites obtained from such choices with $n_{\bullet\bullet}{=}n$ forms a coproduct cocone (the coproduct of coproducts is a coproduct). For a given choice of $r$, $m_i$  and $n_{ij}$ as above one can also form a composite
\[ \xymatrix{{\opE\limits_i\opF\limits_j\opG\limits_k} \ar[r]^-{\opE\limits_i\opc\limits_{jk}} & {\opE\limits_i(\opFoG\limits_{jk})} \ar[r]^-{\opc\limits_{i(jk)}} & {(E{\comp}(F{\comp}G))_{n_{\bullet\bullet}}}} \]
and for all $n \in \N$, the set of all such composites obtained from such choices with $n_{\bullet\bullet}{=}n$ forms a coproduct cocone because $E$ is distributive. Thus for each $n$, $E$, $F$ and $G$ one has a unique isomorphism $\alpha$, such that for all choices of $r$, $m_i$ and $n_{ij}$ with $n_{\bullet\bullet}{=}n$, the diagram
%%%%%
\begin{equation}\label{def-alpha}
\xymatrix{{\opE\limits_i\opF\limits_j\opG\limits_k} \ar[r]^-{\opc\limits_{ij}\opG\limits_k} \ar[dr]_-{\opE\limits_i\opc\limits_{jk}} & {(\opEoF\limits_{ij})\opG\limits_k} \ar[r]^-{\opc\limits_{(ij)k}} & {((E{\comp}F){\comp}G)_n} \ar[d]^{\alpha} \\
& {\opE\limits_i(\opFoG\limits_{jk})} \ar[r]_-{\opc\limits_{i(jk)}} & {(E{\comp}(F{\comp}G))_n}}
\end{equation}
%%%%%
commutes.
\begin{proposition}\label{prop:DistV}
The data $(I,\comp,\alpha,\lambda,\rho)$ just described is a monoidal structure for $\Dist(\ca V)$.
The category $\Mon {\Dist(\ca V)}$ is isomorphic to the category of distributive multitensors and morphisms thereof.
\end{proposition}
\begin{proof}
The case of (\ref{def-alpha}) for which $m_i{=}1$ amounts to the commutativity of the outside of
%%%%%
\[ \xymatrix @R=1.5em @C=1em {& {(\opEoI\limits_i)\opF\limits_k} \ar[rr]^-{\opc\limits_{ik}} && {((E{\comp}I){\comp}F)_n} \ar[dd]^{\alpha} \\ {\opE\limits_i\opF\limits_k} \ar[ur]^{\rho^{-1}\opF\limits_k} \save \POS?="dteq" \restore \ar[rr]^-{\opc\limits_{ik}} \ar[dr]_{\opE\limits_i\lambda^{-1}} \save \POS?="dbeq" \restore && {\opEoF\limits_{ik}} \ar[ur]|{\rho^{-1}{\comp}F} \save \POS?="cteq" \restore \ar[dr]|{E{\comp}\lambda^{-1}} \save \POS?="cbeq" \restore \\ & {\opE\limits_i(\opIoF\limits_k)} \ar[rr]_-{\opc\limits_{ik}} && {(E{\comp}(I{\comp}F))_n}
\POS "dteq"; "cteq" **@{}; ?*{=} \POS "dbeq"; "cbeq" **@{}; ?*{=}} \]
%%%%%
and the inner commutativities indicated here are obtained from the definition of the arrow map of ``$\comp$''. But the
\[ \opc_{ik} : \opE\limits_i\opF\limits_k \rightarrow \opEoF\limits_{ik} \]
for all choices with $n_{\bullet\bullet}{=}n$ form a coproduct cocone, and so the triangle in the above diagram, which is the unit coherence for $\Dist(\ca V)$, must commute also. For $E$, $F$, $G$ and $H$ in $\Dist(\ca V)$ we will now see that the corresponding associativity pentagon commutes. For each $n$ and choice of $r$, $p_i$ for all $1{\leq}i{\leq}r$, $m_{ij}$ for all $i$ and $1{\leq}j{\leq}p_i$, and $n_{ijk}$ for all $i$, $j$ and $1{\leq}k{\leq}m_{ij}$, such that $n_{\bullet\bullet\bullet}{=}n$, we get a diagram of the form:
%%%%%%%%
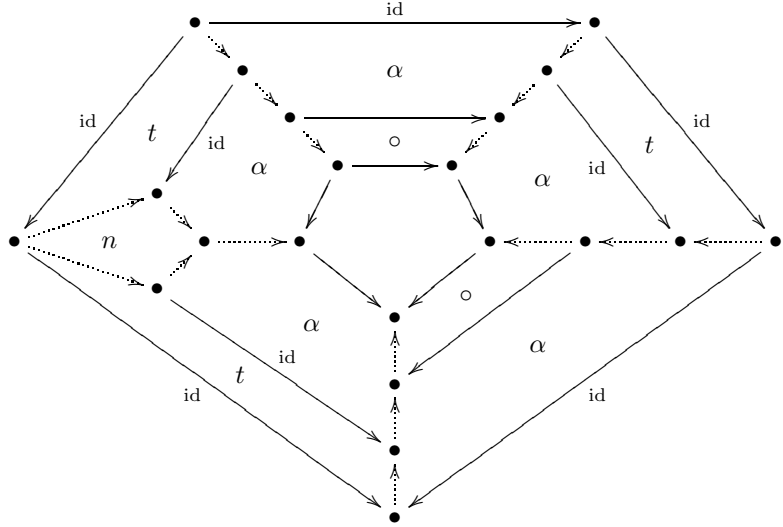
\[ \xygraph{{}="centre" [l(.6)] [u(.8)]{\bullet}="a" "centre" [r(.6)] [u(.8)]{\bullet}="b" "centre" [l]{\bullet}="c" [rr]{\bullet}="d" [d(.8)] [l]{\bullet}="e"
"a":"b"|{}="m1":"d"|{}="m3":"e"|{}="m5" "a":"c"|{}="m2":"e"|{}="m4"
"a" [u(.5)] [l(.5)]{\bullet}="a1" [u(.5)] [l(.5)]{\bullet}="a2" [u(.5)][l(.5)]{\bullet}="a3":@{.>}"a2":@{.>}"a1":@{.>}"a"
"b" [u(.5)] [r(.5)]{\bullet}="b1" [u(.5)][r(.5)]{\bullet}="b2" [u(.5)][r(.5)]{\bullet}="b3":@{.>}"b2":@{.>}"b1":@{.>}"b"
"c" [l]{\bullet}="c1" [ll]{\bullet}="c3" "c1" [u(.5)][l(.5)]{\bullet}="c2u" "c1" [d(.5)][l(.5)]{\bullet}="c2d"
"c3":@{.>}"c2u":@{.>}"c1":@{.>}"c" "c3":@{.>}"c2d":@{.>}"c1"
"d" [r]{\bullet}="d1" [r]{\bullet}="d2" [r]{\bullet}="d3":@{.>}"d2":@{.>}"d1":@{.>}"d"
"e" [d(.7)]{\bullet}="e1" [d(.7)]{\bullet}="e2" [d(.7)]{\bullet}="e3":@{.>}"e2":@{.>}"e1":@{.>}"e"
"a3":"b3"|{}="m12"^{\id}:"d3"|{}="m32"^{\id}:"e3"|{}="m52"^{\id} "a3":"c3"|{}="m22"_{\id}:"e3"|{}="m42"_{\id}
"a2":"c2u"|{}="m21"^{\id} "c2d":"e2"|{}="m41"^{\id} "a1":"b1"|{}="m11" "b2":"d2"|{}="m31"_{\id} "d1":"e1"|{}="m51"
"m1":@{}"m11"|*{\comp}:@{}"m12"|*{\alpha} "m2":@{}"m21"|*{\alpha}:@{}"m22"|*{t}
"m3":@{}"m31"|*{\alpha}:@{}"m32"|*{t} "m4":@{}"m41"|*{\alpha}:@{}"m42"|*{t}
"m5":@{}"m51"|*{\comp}:@{}"m52"|*{\alpha} "c3":@{}"c1"|*{n}} \]
%%%%%%%%
where the inner-most pentagon what we are trying to prove the commutativity of. The outer pentagon has all vertices equal to $\opE\limits_i\opF\limits_j\opG\limits_k\opH\limits_l$. The composites of the dotted paths of length $3$, when taken over all choices, form coproduct cocones of each of the vertices of the inner pentagon. For instance for the top left vertex we have
%%%%%%%%
\[ \xygraph{{\opE\limits_i\opF\limits_j\opG\limits_k\opH\limits_l}
:@{.>}[rr] {(\opEoF\limits_{ij})\opG\limits_k\opH\limits_l}^-{cGH}
:@{.>}[rr] {(\opbEoFboG\limits_{ijk})\opH\limits_l}^-{cH}
:@{.>}[r(2.5)] {(((E{\comp}F){\comp}G){\comp}H)_n}^-{c}} \]
%%%%%%%%
and the two indicated paths involving the left most vertex are
%%%%%%%%
\[ \xygraph{{\opE\limits_i\opF\limits_j\opG\limits_k\opH\limits_l}
:@{.>}[rr] {(\opEoF\limits_{ij})\opG\limits_k\opH\limits_l}^-{cGH}
:@{.>}[r(2.5)] {(\opEoF\limits_{ij})(\opGoH\limits_{kl})}^-{E{\comp}Fc}
:@{.>}[r(2.5)] {(((E{\comp}F){\comp}G){\comp}H)_n}^-{c}} \]
%%%%%%%%
and
%%%%%%%%
\[ \xygraph{{\opE\limits_i\opF\limits_j\opG\limits_k\opH\limits_l}
:@{.>}[rr] {\opE\limits_i\opF\limits_j(\opGoH\limits_{kl})}^-{EFc}
:@{.>}[r(2.5)] {(\opEoF\limits_{ij})(\opGoH\limits_{kl})}^-{cG{\comp}H}
:@{.>}[r(2.5)] {(((E{\comp}F){\comp}G){\comp}H)_n}^-{c}} \]
%%%%%%%%
and in a similar vein the reader will easily supply the details of the other dotted paths. The labels of the regions of the diagram indicate why the corresponding region commutes: ``$\alpha$'' means the region commutes by the definition of $\alpha$, ``n'' indicates commutativity because of naturality, ``$\comp$'' indicates commutativity because of the definition of the arrow map of $\comp$, and ``t'' indicates that the region commutes trivially. The outer pentagon of course also commutes trivially. Since all this is true for all choices of the $r$, $p_i$, $m_{ij}$ and $n_{ijk}$, we obtain the commutativity of the inner pentagon since the top left dotted composites together exhibit $(((E{\comp}F){\comp}G){\comp}H)_n$ as a coproduct. The statement about $\Mon {\Dist(\ca V)}$ follows immediately by unpacking the definitions involved.
\end{proof}
%

%%%%%%%%%%%%%%%%%%%%%%%%%%%%%%%%%%%%%%%%%%%%%%%%%%
%%%%%%%%%%%%%%%%%%%%%%%%%%%%%%%%%%%%%%%%%%%%%%%%%%
\section{Monads from multitensors}\label{ssec:mult-fmon}

Multitensors generalise non-symmetric operads by example(\ref{ex:nso-as-mt}). Given certain hypotheses on the ambient braided monoidal category $\ca V$, a non-symmetric operad therein gives rise to a monad on $\ca V$ whose algebras are those of the original operad. Thus one is led to ask whether one can define a monad from a multitensor in a similar way. Such a construction is described in the present section, and we continue to assume throughout this section that $\ca V$ has coproducts.

Define the functor $\Gamma : \Dist(\ca V) \rightarrow \End(\ca V)$ as
\[ \Gamma(E)(X) = \coprod_{n{\geq}0} \opE_{1{\leq}i{\leq}n}X \]
and so for each $X$ in $\ca V$ we get
\[ c_n : \opE_{1{\leq}i{\leq}n}X \rightarrow \Gamma(E)(X) \]
for $n \in \N$ forming a coproduct cocone. By the definition of $I$ the map $c_1:X{\rightarrow}\Gamma(I)(X)$ is an isomorphism, and we define that the inverses of these maps are the components of an isomorphism $\gamma_0 : 1_{\ca V}{\rightarrow}\Gamma(I)$. For $X$ in $\ca V$ and $m$ and $n_i$ in $\N$ where $1{\leq}i{\leq}m$, we can consider composites
\[ \xygraph{!{(0,0);(2,0):(0,1)::}
{\opE\limits_i\opF\limits_jX}:[r]{\opE\limits_i\Gamma{F}X}^-{\opE\limits_i\opc\limits_j}:[r]{\Gamma(E)\Gamma(F)X}^-{c_m}} \]
and since $E$ is distributive all such composites exhibit $\Gamma(E)\Gamma(F)X$ as a coproduct. For $X$, $m$ and $n_i$ as above one also has composites
\[ \xygraph{!{(0,0);(2,0):(0,1)::}
{\opE\limits_i\opF\limits_jX}:[r]{(\opEoF\limits_{ij})X}^-{\opc\limits_{ij}}:[r]{\Gamma(E{\comp}F)X}^-{c_{n_{\bullet}}}} \]
and all such composites exhibit $\Gamma(E{\comp}F)X$ as a coproduct. Thus there is a unique isomorphism $\gamma_2$ making
%%%%%%%%
\[ \xygraph{!{(0,0);(0,1.5):(0,1.5)::}
{\opE\limits_i\opF\limits_jX}="a" ([d]{(\opEoF\limits_{ij})X}="c" [d][r(.5)]{\Gamma(E{\comp}F)X}="e",
[r]{\opE\limits_i\Gamma(F)X}="b" [d]{\Gamma(E)\Gamma(F)X}="d")
"a":"b"^-{\opE\limits_i\opc\limits_j}:"d"^-{c_m}:"e"^-{\gamma_2}
"a":"c"_-{\opc\limits_{ij}}:"e"_-{c_{n_\bullet}}} \]
%%%%%%%%
commute, and $\gamma_2$ is clearly natural in $X$.
\begin{proposition}\label{prop:Gamma}
The data $(\gamma_0,\gamma_2)$ make $\Gamma$ into a monoidal functor. For any distributive multitensor $E$, one has an isomorphism $\Mon E \iso \Alg {\Gamma E}$ commuting with the forgetful functors into $\ca V$.
\end{proposition}
\begin{proof}
The definition of $\gamma_2$ in the case where $E{=}I$ and the $m{=}1$ says that the outside of
%%%%%%%%
\[ \xygraph{!{(0,0);(0,1.5):(0,1.5)::}
{\opF\limits_jX}="a" ([d]{(\opIoF\limits_j)X}="c" [d][r(.5)]{\Gamma(I{\comp}F)X}="e",
[r]{\Gamma(F)X}="b" [d]{\Gamma(I)\Gamma(F)X}="d")
"a":"b"^-{\opc\limits_j}:"d"^-{\gamma_0\Gamma{F}}:"e"^-{\gamma_2}
"a":"c"_-{\lambda^{-1}}:"e"_-{c_n}
"b":@/_{1.5pc}/"e"^{\Gamma\lambda^{-1}}="m" "a"-@{}"m"|(.6)*{=}} \]
%%%%%%%%
commutes for all $m \in \N$, and the region labelled with ``$=$'' commutes because of the definition of the arrow maps of $\comp$. Thus the inner triangle, which is the left unit monoidal functor coherence axiom, commutes also. The definition of $\gamma_2$ in the case where $F{=}I$ and the $n_i$'s are all $1$ says that the outside of
%%%%%%%%
\[ \xygraph{!{(0,0);(0,2):(0,1.5)::}
{\opE\limits_i\opI\limits_jX}="a" ([d]{(\opEoI\limits_{ij})X}="c" [d][r(.5)]{\Gamma(E{\comp}I)X}="e",
[r]{\opE\limits_i\Gamma(I)X}="b" [d]{\Gamma(E)\Gamma(I)X}="d")
"a":"b"^-{\opE\limits_i\gamma_0}:"d"^-{c_m}:"e"^-{\gamma_2}
"a":"c"_-{\rho^{-1}}:"e"_-{c_{n_\bullet}}
"c"-@{}"d"|(.4){}="mp" "mp"[u(.4)]{\Gamma(E)X}="m"
"a":"m"^-{c_m}(:"e"^{\Gamma\rho^{-1}}|(.3){}="va",:"d"^-{\Gamma(E)\gamma_0})
"c"-@{}"va"|(.6)*{=} "m"-@{}"b"|(.45)*{=}} \]
%%%%%%%%
commutes for all $m \in \N$, and the regions labelled with ``$=$'' commute because of the definition of the arrow maps of $\comp$. Thus the inner triangle, which is the right unit monoidal functor coherence axiom, commutes also. So it remains to verify that for $E$, $F$ and $G$ in $\Dist(\ca V)$, that
%%%%%%%%
\begin{equation}\label{eq:pent}
\xygraph{!{0;(2.5,0):(0,.6)::}
{} ([l(.6)] [u(.8)]{\Gamma(E)\Gamma(F)\Gamma(G)}="a", [r(.6)] [u(.8)]{\Gamma(E{\comp}F)\Gamma(G)}="b", [l]{\Gamma(E)\Gamma(F{\comp}G)}="c" [rr]{\Gamma((E{\comp}F){\comp}G)}="d"
[d(.8)] [l]{\Gamma(E{\comp}(F{\comp}G))}="e") "a":"b"^{\gamma_2\Gamma(G)}:"d"^{\gamma_2}:"e"^{\Gamma\alpha} "a":"c"_{\Gamma(E)\gamma_2}:"e"_{\gamma_2}}
\end{equation}
%%%%%%%%
commutes. Now given $X$ in $\ca V$ and $r$, $m_i$ and $n_{ij}$ in $\N$ where $1{\leq}i{\leq}r$ and $1{\leq}j{\leq}m_i$, one obtains a diagram of the form
%%%%%%%%
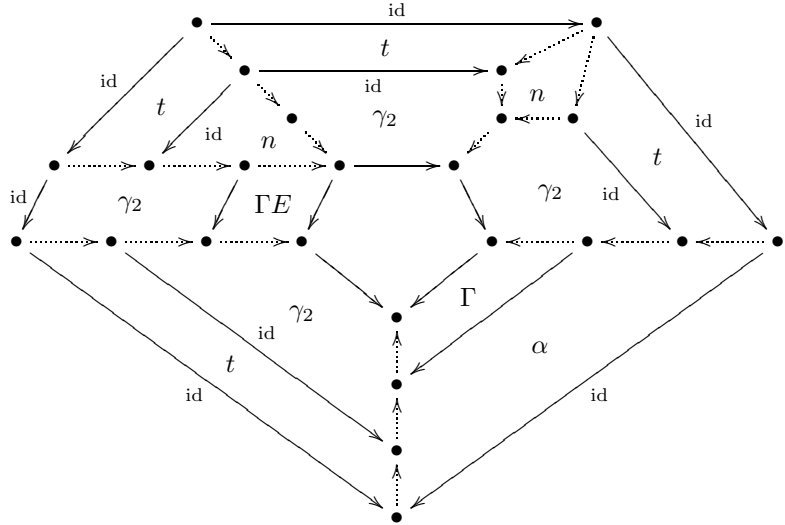
\[ \xygraph{{}="centre" [l(.6)] [u(.8)]{\bullet}="a" "centre" [r(.6)] [u(.8)]{\bullet}="b" "centre" [l]{\bullet}="c" [rr]{\bullet}="d" [d(.8)] [l]{\bullet}="e"
"a":"b"|{}="t":"d"|{}="tr":"e"|{}="br" "a":"c"|{}="tl":"e"|{}="bl"
"a" [u(.5)] [l(.5)]{\bullet}="a1" [u(.5)] [l(.5)]{\bullet}="a2" [u(.5)][l(.5)]{\bullet}="a3":@{.>}"a2":@{.>}"a1":@{.>}"a"
"a" [l]{\bullet}="a1l"  [l]{\bullet}="a2l" [l]{\bullet}="a3l":@{.>}"a2l":@{.>}"a1l":@{.>}"a"
"b" [u(.5)] [r(.5)]{\bullet}="b1" ([ur]{\bullet}="b3",[u(.5)]{\bullet}="b2l", [r(.75)]{\bullet}="b2r")
"b3":@{.>}"b2l":@{.>}"b1":@{.>}"b" "b3":@{.>}"b2r":@{.>}"b1"
"c" [l]{\bullet}="c1" [l]{\bullet}="c2" [l]{\bullet}="c3":@{.>}"c2":@{.>}"c1":@{.>}"c"
"d" [r]{\bullet}="d1" [r]{\bullet}="d2" [r]{\bullet}="d3":@{.>}"d2":@{.>}"d1":@{.>}"d"
"e" [d(.7)]{\bullet}="e1" [d(.7)]{\bullet}="e2" [d(.7)]{\bullet}="e3":@{.>}"e2":@{.>}"e1":@{.>}"e"
"a3":"b3"^{\id}|{}="t3":"d3"^{\id}|{}="tr3":"e3"^{\id}|{}="br3"
"a3":"a3l"_{\id}|{}="tl3u":"c3"_{\id}|{}="tl3l":"e3"_{\id}|{}="bl3"
"a2":"b2l"_{\id}|{}="t2" "b2r":"d2"_{\id}|{}="tr2" "a2":"a2l"^{\id}|{}="tl2u" "c2":"e2"^{\id}|{}="bl2" "d1":"e1"|{}="br1" "a1l":"c1"|{}="tl1l"
"t3"-@{}"t2"|*{t}-@{}"t"|*{\gamma_2} "b2l"-@{}"b2r"|*{n} "tr3"-@{}"tr2"|*{t}-@{}"tr"|*{\gamma_2}
"br3"-@{}"br1"|*{\alpha}-@{}"br"|*{\Gamma} "tl3u"-@{}"tl2u"|*{t}-@{}"a"|*{n}
"tl3l"-@{}"tl1l"|*{\gamma_2}-@{}"tl"|*{\Gamma{E}} "bl3"-@{}"bl2"|*{t}-@{}"bl"|*{\gamma_2}} \]
%%%%%%%%
where the inner-most pentagon is (\ref{eq:pent}) instantiated at $X$, and all the outer vertices are $\opE\limits_i\opF\limits_j\opG\limits_kX$. The two 3-fold paths into $\Gamma(E)\Gamma(F)\Gamma(G)(X)$ are the top-leftmost path
\[ \xygraph{{\opE\limits_i\opF\limits_j\opG\limits_kX}:@{.>}[rr]{\opE\limits_i\opF\limits_j\Gamma(G)X}^-{\opE\limits_i\opF\limits_j\opc\limits_k}:@{.>}[r(2)]{\Gamma(E)\opF\limits_j\Gamma(G)X}^-{\opc\limits_i}:@{.>}[rrr]{\Gamma(E)\Gamma(F)\Gamma(G)(X)}^-{\Gamma(E)\opc_j}} \]
and
\[ \xygraph{{\opE\limits_i\opF\limits_j\opG\limits_kX}:@{.>}[r(2)]{\opE\limits_i\opF\limits_j\Gamma(G)X}^-{\opE\limits_i\opF\limits_j\opc\limits_k}:@{.>}[r(2.5)]{\opE\limits_i\Gamma(F)\Gamma(G)(X)}^-{\opE\limits_i\opc\limits_j}:@{.>}[r(2.8)]{\Gamma(E)\Gamma(F)\Gamma(G)(X)}^-{\opc\limits_i}} \]
and these are equal because of naturality. The composites so formed by taking all choices of $r$, $m_i$ and $n_{ij}$ exhibit $\Gamma(E)\Gamma(F)\Gamma(G)(X)$ as a coproduct because $E$ and $F$ are distributive. The left-most dotted path into $\Gamma(E{\comp}F)\Gamma(G)(X)$ is
\[ \xygraph{{\opE\limits_i\opF\limits_j\opG\limits_kX}:@{.>}[r(2)]{\opE\limits_i\opF\limits_j\Gamma(G)X}^-{\opE\limits_i\opF\limits_j\opc\limits_k}:@{.>}[r(2)]{(\opEoF\limits_{ij})\Gamma(G)X}^-{\opc\limits_{ij}}:@{.>}[r(2.5)]{\Gamma(E{\comp}F)\Gamma(G)(X)}^-{\opc\limits_{ij}}}, \]
the other path into $\Gamma(E{\comp}F)\Gamma(G)(X)$ is
\[ \xygraph{{\opE\limits_i\opF\limits_j\opG\limits_kX}:@{.>}[r(1.8)]{\opEoF\limits_{ij}\opG\limits_kX}^-{\opc\limits_{ij}}:@{.>}[r(2.5)]{(\opEoF\limits_{ij})\Gamma(G)X}^-{\opEoF\limits_{ij}\opc\limits_k}:@{.>}[r(2.5)]{\Gamma(E{\comp}F)\Gamma(G)(X)}^-{\opc\limits_{ij}}}, \]
and similarly the reader will easily supply the definitions of the other dotted paths in the above diagram. The labelled regions of that diagram commute for the reasons indicated by the labels as with the proof of proposition(\ref{prop:DistV}), the region labelled by ``$\Gamma$'' commutes by the definition of the arrow map of $\Gamma$, and the region labelled by ``$\Gamma{E}$'' commutes by the definition of the arrow map of $\Gamma{E}$. The outer diagram commutes trivially and since this is all true for all choices of the $r$, $m_i$ and $n_{ij}$, the inner pentagon commutes as required. The statement about $\Mon E$ follows immediately by unpacking the definitions involved.
\end{proof}
\begin{example}\label{ex:nso-as-mtII}
One can apply proposition(\ref{prop:Gamma}) to the case of example(\ref{ex:nso-as-mt}) when $(\ca V,\tensor,I)$ is a distributive braided monoidal category, because then the multitensor on $\ca V$ determined by a non-symmetric operad will also be distributive. In this way one obtains the usual construction of the monad induced by a non-symmetric operad.
\end{example}
\begin{example}\label{ex:Tcross-dist-II}
Applying proposition(\ref{prop:Gamma}) to the case of a distributive monoidal category $(\ca V,\tensor,I)$ as in
example(\ref{ex:Tcross-dist}), one recovers the usual monoid monad $M{:=}\Gamma(\tensor)$. In the case where $\tensor$ is cartesian product and $T$ preserves coproducts, in view of $\Gamma(T^{\times}){=}MT$ one obtains a monad structure on $MT$, and thus a monad distributive law $\lambda:TM{\rightarrow}MT$, and the algebras of $MT$ are monoids in $\Alg T$ by proposition(\ref{prop:Tcross-cats}). In terms of $\Gamma$ and $T^{\times}$ one can describe $\lambda$ explicitly. The substitution for $T^{\times}$, described in example(\ref{ex:mult-Tcross}), is a map $\mu^{\times}:T^{\times}{\comp}T^{\times}{\rightarrow}T^{\times}$ in $\Dist(\ca V)$, and $\lambda$ is the composite
\[ \xymatrix{{TM} \ar[r]^-{{\eta}TM\eta} & {MTMT} \ar[r]^-{\Gamma\mu^{\times}} & {MT}} \]
in $\End(\ca V)$.
\end{example}
%

%%%%%%%%%%%%%%%%%%%%%%%%%%%%%%%%%%%%%%%%%%%%%%%%%%
%%%%%%%%%%%%%%%%%%%%%%%%%%%%%%%%%%%%%%%%%%%%%%%%%%
\section{Multitensors as operads}\label{ssec:masop}

Given a cartesian monad $T$ on a finitely complete category $\ca V$ one has the well-known notion of \emph{$T$-operad} as described for example in \cite{Lei}. There is an analogous notion of $T$-multitensor and we shall describe this in the present section. Under certain conditions the given monad $T$ distributes with the monoid monad $M$ on $\ca V$ and the composite monad $MT$ is again cartesian, in which case one has an equivalence of categories between $T$-multitensors and $MT$-operads. The theory described in this section requires that $T$ is a little more than cartesian, namely that it is p.r.a in the sense of \cite{Fam2fun}, and that $\ca V$ is lextensive. Both notions will be recalled here for the readers' convenience.

We recall some aspects of the theory parametric right adjoints from \cite{Fam2fun}. A functor $T:{\ca A}{\rightarrow}{\ca B}$ is a \emph{parametric right adjoint} (p.r.a){\footnotemark{\footnotetext{We reserve the right to use this abbreviation also as an adjective, as in ``$T$ is {\bf p}arametrically {\bf r}epresent{\bf a}ble''.}}} when for all $A \in \ca A$, the induced functors
\[ T_A : \ca A/A \rightarrow \ca B/TA \]
given by applying $T$ to arrows have left adjoints, and when $\ca A$ has a terminal object $1$, this is equivalent to asking that $T_1$ has a left adjoint. Right adjoints are clearly p.r.a and p.r.a functors are closed under composition. Moreover one has the following simple observation which we shall use often in this work.
\begin{lemma}\label{lem:pra-prod}
Let $I$ be a set and $F_i:\ca A_i{\rightarrow}\ca B_i$ for $i \in I$ be a family of p.r.a functors. Then
\[ \xymatrix{{\prod\limits_i\ca A_i} \ar[r]^-{\prod\limits_iF_i} & {\prod\limits_i\ca B_i}} \]
is p.r.a.
\end{lemma}
\begin{proof}
Given $X_i \in \ca A_i$ for $i \in I$, we have $(\prod\limits_iF_i)_{(X_i)}{=}\prod\limits_i((F_i)_{X_i})$, which as a product of right adjoints is a right adjoint.
\end{proof}
\noindent There is a more explicit characterisation of p.r.a functors which is sometimes useful. A map $f:B{\rightarrow}TA$ is \emph{$T$-generic} when for any $\alpha$, $\beta$, and $\gamma$ making the outside of
\[ \xymatrix{B \ar[r]^-{\alpha} \ar[d]_{f} & {TX} \ar[d]^{T\gamma} \\
{TA} \ar[r]_-{T\beta} \ar@{.>}[ur]|{T\delta} & {TZ}} \]
commute, there is a unique $\delta$ for which $\gamma \comp \delta = \beta$ and $T({\delta}) \comp f = \alpha$. The alternative characterisation says that $T$ is p.r.a iff every map $f:B{\rightarrow}TA$ factors as
\[ \xymatrix{B \ar[r]^-{g} & {TC} \ar[r]^-{Th} & {TA}} \]
where $g$ is generic, and such generic factorisations are unique up to isomorphism if they exist (see \cite{Fam2fun} for more details). One defines a monad $(T,\eta,\mu)$ on a category $\ca V$ to be p.r.a when $T$ is p.r.a as a functor, and $\eta$ and $\mu$ are cartesian transformations. One has the following corresponding definition for multitensors.
\begin{definition}
A multitensor $(E,u,\sigma)$ on $\ca V$ is p.r.a when $E:\ca M\ca V{\rightarrow}\ca V$ is p.r.a and $u$ and $\sigma$ are cartesian transformations.
\end{definition}
\noindent It is straight-forward to observe that $E$ is p.r.a iff $E_n:\ca V^n{\rightarrow}\ca V$ is p.r.a for each $n \in \N$.
\begin{example}\label{ex:Tcross-pra}
Let $(T,\eta,\mu)$ be a p.r.a monad on $\ca V$ a category with finite products. First note that $T^{\times}_n$ is the composite
\[ \xymatrix{{\ca V^n} \ar[r]^-{T^n} & {\ca V^n} \ar[r]^-{\prod} & {\ca V}} \]
and so is p.r.a. by lemma(\ref{lem:pra-prod}) and the composability of p.r.a's. From \cite{Fam2fun} lemma(2.14) the canonical maps 
\[ \begin{array}{c} {k_{X_i} : T \prod\limits_iX_i \rightarrow \prod\limits_iTX_i} \end{array} \]
which measure the extent to which $T$ preserves products are cartesian natural in the $X_i$. Thus $T^{\times}$ is a p.r.a multitensor.
\end{example}
\noindent For a p.r.a monad $(T,\eta,\mu)$ on a category $\ca V$ recall that a \emph{$T$-operad} is cartesian monad morphism $\alpha:A{\rightarrow}T$. That is, $A$ is a monad on $\ca V$, $\alpha$ is a natural transformation $A{\rightarrow}T$ which is compatible with the monad structures, and the naturality squares of $\alpha$ are pullbacks. The cartesianness of $\alpha$ and p.r.a'ness of $T$ implies that $A$ is itself a p.r.a monad. For instance when $T{=}\ca T$ the monad on the category $\PSh {\G}$ of globular sets whose algebras are strict $\omega$-categories, to be recalled in detail in section(\ref{sec:wcat-monad}), $T$-operads are the $\omega$-operads of Batanin \cite{Bat98}. By analogy one has the following definition for multitensors.
\begin{definition}
Let $(T,\eta,\mu)$ be a p.r.a monad on $\ca V$ a category with finite products. A \emph{$T$-multitensor} is a cartesian multitensor morphism $\varepsilon:E{\rightarrow}T^{\times}$.
\end{definition}
\begin{example}\label{ex:0mult}
We will now unpack this notion in the case where $\ca V=\Set$ and $T$ is the identity monad. Because of the pullback squares
%%%%%%%%
\[ \xygraph{!{0;(1.5,0):(0,.8)::}
{\opE\limits_iX_i}="tl" [r] {\prod\limits_iTX_i}="tr" [d] {(T1)^n}="br" [l] {E_n1}="bl"
"tl" (:"tr"^{\varepsilon_{X_i}}:"br"^{\prod\limits_iTt_{X_i}},
:"bl"_{\opE\limits_it_{X_i}}:"br"_{\varepsilon_1})
"tl" [d(.3)r(.3)] (:@{-}[l(.15)],:@{-}[u(.15)])} \]
%%%%%%%%
the data for $\varepsilon$ amounts to a sequence of objects $\overline{E}_n:=E_n1 \in \ca V$ for $n \in \N$, together with maps $\varepsilon_{n,i}:\overline{E}_n{\rightarrow}T1$ for $1{\leq}i{\leq}n$. In this case $T1{=}1$ so $\varepsilon$ amounts to a sequence $(\overline{E}_n:n \in \N)$ of sets. In terms of this data one has
\begin{equation}\label{eq:op-deform}
\begin{array}{c} {\opE\limits_{1{\leq}i{\leq}n}X_i = \overline{E}_n \times \prod\limits_iX_i} \end{array}
\end{equation}
The unit of the multitensor amounts to an element $\overline{u}:1{\rightarrow}\overline{E}_1$, and the substitution $\sigma$ amounts to functions
\[ \overline{\sigma}_{n_1,...,n_k} : \overline{E}_k \times \overline{E}_{n_1} \times ... \times \overline{E}_{n_k} \rightarrow \overline{E}_{n_{\bullet}} \]
for each finite sequence $(n_1,...,n_k)$ of natural numbers. The multitensor axioms for $(E,u,\sigma)$ correspond to axioms that make $(\overline{E},\overline{u},\overline{\sigma})$ a non-symmetric operad in $\Set$.
\end{example}
We assume throughout this section that $\ca V$ is lextensive. Let us now recall this notion. A category $\ca V$ is \emph{lextensive}{\footnotemark{\footnotetext{Usually lextensivity is defined using only finite coproducts whereas we work with small ones.}}} \cite{CLW} when it has finite limits, coproducts and for each family of objects $(X_i : i \in I)$ of $\ca V$ the functor
\[ \begin{array}{c} {\prod\limits_{i{\in}I} \ca V/X_i \rightarrow \ca V/{\left(\coprod\limits_{i{\in}I} X_i\right)}}
\end{array} \]
which sends a family of maps $(h_i:Z_i{\rightarrow}X_i)$ to their coproduct is an equivalence. This last property is equivalent to saying that $\ca V$ has a strict initial object and that coproducts in $\ca V$ are disjoint and stable. There are many examples of lextensive categories: for instance every Grothendieck topos is lextensive, as is $\CAT$. Moreover if $T$ is a coproduct preserving monad on a lextensive category $\ca V$ then $\Alg T$ is also lextensive: for such a $T$ the forgetful functor $\Alg T{\rightarrow}\ca V$ creates finite limits and coproducts, and so these exist in $\Alg T$ and interact as nicely as they did in $\ca V$. Thus in particular the category of algebras of any higher operad is lextensive. Note in particular that lextensivity implies distributivity (see \cite{CLW}) and so the results of the previous two sections apply in this one. The next result summarises how lextensivity interacts well with p.r.a'ness.
\begin{lemma}\label{lem:lext-pra}
Let $\ca A$ and $\ca B$ be lextensive and $I$ be a set.
\begin{enumerate}
\item  The functor $\coprod:\ca A^I{\rightarrow}\ca A$, which takes an $I$-indexed family of objects of $\ca A$ to its coproduct, is p.r.a.\label{lpra1}
\item  If $F_i:\ca A{\rightarrow}\ca B$ for $i \in I$ are p.r.a functors, then $\coprod\limits_iF_i:\ca A{\rightarrow}\ca B$ is p.r.a.\label{lpra2}
\item  If $F_i:\ca A{\rightarrow}\ca B$ for $i \in I$ are functors and $\phi_i:F_i{\rightarrow}G_i$ are cartesian transformations, then $\coprod\limits_i\phi_i:\coprod\limits_iF_i{\rightarrow}\coprod\limits_iG_i$ is cartesian.\label{lpra3}
\end{enumerate}
\end{lemma}
\begin{proof}
(\ref{lpra1}): given a family $(X_i:i{\in}I)$ of objects of $\ca A$, the functor $(\coprod)_{(X_i)}$ is just the functor
\[ \begin{array}{c} {\prod\limits_{i{\in}I} \ca A/X_i \rightarrow \ca A/{\left(\coprod\limits_{i{\in}I} X_i\right)}}
\end{array} \]
which is an equivalence, and thus a right adjoint.\\
(\ref{lpra2}): $\coprod\limits_iF_i$ is the composite
\[ \xymatrix{{\ca A} \ar[r]^-{\Delta} & {\ca A^I} \ar[r]^-{\prod\limits_iF_i} & {\ca B^I} \ar[r]^-{\coprod} & {\ca B}} \]
of a right adjoint (since $\ca A$ has coproducts) followed by a p.r.a (by lemma(\ref{lem:pra-prod}) followed by another p.r.a (by (\ref{lpra1}), and so is p.r.a.\\
(\ref{lpra3}): the naturality square for $\coprod\limits_i\phi_i$ corresponding to $f:X{\rightarrow}Y$ in $\ca A$ is the coproduct of the cartesian naturality squares
\[ \xymatrix{{F_iX} \ar[r]^-{\phi_{i,X}} \ar[d]_{F_if} & {G_iX} \ar[d]^{G_if} \\ {F_iY} \ar[r]_-{\phi_{i,Y}} & {G_iY}} \] 
and so by (\ref{lpra1}) is itself a pullback.
\end{proof}
\noindent Denote by $\PraDist(\ca V)$ and $\PraEnd(\ca V)$ the subcategories of $\Dist(\ca V)$ and $\End(\ca V)$ respectively, whose objects are p.r.a's and arrows are cartesian transformations.
\begin{proposition}\label{prop:pra-restrict}
Let $\ca V$ be lextensive. The monoidal structure of $\Dist(\ca V)$ restricts to $\PraDist(\ca V)$, and $\Gamma$ restricts to a strong monoidal functor \[ \PraDist(\ca V){\rightarrow}\PraEnd(\ca V) \] (which we shall also denote by $\Gamma$).
\end{proposition}
\begin{proof}
Any functor $1{\rightarrow}\ca A$ out of the terminal category is p.r.a, and thus one readily verifies that the functors $\ca V^n{\rightarrow}\ca V$ constant at the initial object $0$ of $\ca V$ are p.r.a also. Since $1_{\ca V}$ is p.r.a the unit of $\Dist(\ca V)$ is p.r.a. For p.r.a $E$ and $F \in \Dist(\ca V)$ we must verify that $E{\comp}F$ is p.r.a. By the formula
\[ (E \comp F)_n = \coprod_{n_1+...+n_k{=}n} E_k(F_{n_1},...,F_{n_k}) \]
and lemma(\ref{lem:lext-pra}) it suffices to show that each summand is p.r.a. But $E_k(F_{n_1},...,F_{n_k})$ is the composite
\[ \xymatrix{{\prod\limits_i\ca V^{n_i}} \ar[r]^-{\prod\limits_iF_{n_i}} & {\ca V^k} \ar[r]^-{E_k} & {\ca V}} \]
which is p.r.a by lemma(\ref{lem:pra-prod}). Given $\varepsilon:E{\rightarrow}E'$ and $\phi:F{\rightarrow}F'$ in $\PraDist(\ca V)$ we must show that $\varepsilon{\comp}\phi$ is cartesian. By lemma(\ref{lem:lext-pra}) it suffices to show that
\[ \xymatrix{{E_k(F_{n_1},...,F_{n_k})} \ar[rr]^-{\varepsilon_k(\phi_{n_1},...,\phi_{n_k})} && {E'_k(F'_{n_1},...,F'_{n_k})}} \]
is cartesian. But this natural transformation is the composite
%%%%%%%%
\[ \xygraph{!{0;(2,0):(0,.5)::}
{\prod\limits_i\ca V^{n_i}}="a" [r] {\ca V^k}="b" [r] {\ca V}="c"
"a":@/^{1.5pc}/"b"^-{\prod\limits_iF_{n_i}}|(.4){}="l"
"a":@/_{1.5pc}/"b"_-{\prod\limits_iF'_{n_i}}
"b":@/^{1.5pc}/"c"^-{E_k}|{}="r" "b":@/_{1.5pc}/"c"_-{E'_k}
"l" [d(.3)] :@{=>}^{\prod\limits_i\phi_{n_i}} [d(.4)] "r" [d(.3)] :@{=>}^{\varepsilon_k} [d(.4)]} \]
%%%%%%%%
and so as a horizontal composite of cartesian transformations between pullback preserving functors, is indeed cartesian. Thus the monoidal structure of $\Dist(\ca V)$ restricts to $\PraDist(\ca V)$, and to finish the proof we must verify that $\Gamma$ preserves p.r.a objects and cartesian transformations. Let $E \in \Dist(\ca V)$ be p.r.a. By lemma(\ref{lem:lext-pra}), to establish that $\Gamma(E)$ is p.r.a it suffices to show that for all $n \in \N$, the functor $X \mapsto E_n(X,...,X)$ is p.r.a, but this is just the composite
\[ \xymatrix{{\ca V} \ar[r]^-{\Delta} & {\ca V^n} \ar[r]^-{E_n} & {\ca V}} \]
which is p.r.a since $E_n$ is. Let $\phi:E{\rightarrow}F$ in $\Dist(\ca V)$ be cartesian and let us see that $\Gamma(\phi)$ is cartesian. By lemma(\ref{lem:lext-pra}) this comes down to the cartesian naturality in $X$ of the maps
\[ \phi_{n,X,...,X} : E_n(X,...,X) \rightarrow F_n(X,...,X) \]
which is an instance of the cartesianness of $\phi_n$.
\end{proof}
\begin{example}\label{ex:Tcross-dist-III}
From examples(\ref{ex:Tcross-pra}) and example(\ref{ex:Tcross-dist}) $T^{\times}$ is a p.r.a distributive multitensor when $T$ is a coproduct preserving p.r.a monad on a lextensive category $\ca V$. By proposition(\ref{prop:pra-restrict}), the monad $MT$ described in example(\ref{ex:Tcross-dist-II}) is p.r.a and the distributive law $\lambda:TM{\rightarrow}MT$ is cartesian.
\end{example}
Modulo one last digression we are now ready to exhibit the equivalence between $T$-multitensors and $MT$-operads as promised at the beginning of this section. Recall that if $\ca W$ is a monoidal category and $(M,i,m)$ a monoid therein, that the slice $\ca W/M$ gets a canonical monoidal structure. The unit is the unit $i:I{\rightarrow}M$ of the monoid, the tensor product of arrows $\alpha:A{\rightarrow}M$ and $\beta:B{\rightarrow}M$ is the composite
\[ \xymatrix{{A{\tensor}B} \ar[r]^-{\alpha{\tensor}\beta} & {M{\tensor}M} \ar[r]^-{m} & M} \]
and the coherences are inherited from $\ca W$ so that the forgetful functor $\ca W/M{\rightarrow}\ca W$ is strict monoidal. To give $\alpha:A{\rightarrow}M$ a monoid structure in $\ca W/M$ is the same as giving $A$ a monoid structure for which $\alpha$ becomes a monoid homomorphism, and this is just the object part of an isomorphism $\Mon {\ca W/M}{\iso}\Mon {\ca W}/M$ commuting with the forgetful functors into $\ca W$. Moreover given a monoidal functor $F:{\ca W}{\rightarrow}{\ca W'}$, $FM$ is canonically a monoid and one has a commutative square
\[ \xymatrix{{\ca W/M} \ar[r]^-{F_M} \ar[d] & {\ca W'/FM} \ar[d] \\ {\ca W} \ar[r]_-{F} & {\ca W'}} \]
of monoidal functors.

Applying these observations to $\Gamma:\PraDist(\ca V){\rightarrow}\PraEnd(\ca V)$ one obtains for each p.r.a distributive multitensor $E$, a monoidal functor
\[ \Gamma_E : \PraDist(\ca V)/E \rightarrow \PraEnd(\ca V)/\Gamma{E}. \]
An object of $\PraDist(\ca V)/E$ amounts to a functor $A:\ca M\ca V{\rightarrow}\ca V$ together with a cartesian transformation $\alpha:A{\rightarrow}E$. Given such data the distributivity of $A$ is a consequence of the cartesianness of $\alpha$, the distributivity of $E$ and the stability of $\ca V$'s coproducts. The p.r.a'ness of $A$ is also a consequence, because the domain of any cartesian transformation into a p.r.a functor is again p.r.a. A morphism in $\PraDist(\ca V)/E$ from $\alpha$ to $\beta:B{\rightarrow}E$ is just a natural transformation $\phi:A{\rightarrow}B$ such that $\beta\phi{=}\alpha$, because by the elementary properties of pullbacks $\phi$ is automatically cartesian. Thus a monoid in $\PraDist(\ca V)/E$ is simply a cartesian multitensor morphism into $E$. Similarly a monoid in $\PraEnd(\ca V)/\Gamma{E}$ is just a cartesian monad morphism into $\Gamma{E}$, and so by observing its effect on monoids in the case $E{=}T^{\times}$ where $T$ is a coproduct pres!
 
 erving p.r.a monad on $\ca V$, one has a functor
\[ \Gamma_T : \Mult T \rightarrow \Op {MT} \]
from the category of $T$-multitensors to the category of $MT$-operads.
\begin{theorem}\label{thm:TMult-MTOp}
Let $\ca V$ be lextensive and $T$ a coproduct preserving p.r.a monad on $\ca V$. Then the functor $\Gamma_T$ just described is an equivalence of categories $\Mult T \catequiv \Op {MT}$.
\end{theorem}
\begin{proof}
By the way we have set things up it suffices to show that for any p.r.a distributive multitensor $E$ on $\ca V$, the functor $\Gamma_E:{\PraDist(\ca V)/E}{\rightarrow}{\PraEnd(\ca V)/\Gamma{E}}$ is essentially surjective on objects and fully faithful. Let $\alpha:A{\rightarrow}\Gamma(E)$ be a cartesian transformation. Choosing pullbacks
%%%%%%%%
\[ \xygraph{!{0;(2,0):(0,.6)::}
{\opA\limits_iX_i}="abar" [d] {\opE\limits_iX_i}="e" [r] {E_n(1,...,1)}="e1"
[r] {\Gamma{E}(1)}="ge1" [u] {A1}="a1"
"abar" (:"e"_-{\overline{\alpha}_{X_i}} :"e1"_-{\opE\limits_it_{X_i}} :"ge1"_-{c_n}, :"a1" :"ge1"^{c_n},
[d(.3)r(.15)] :@{-} [r(.15)] :@{-} [u(.15)])} \]
%%%%%%%%
for each finite sequence $(X_i:1{\leq}i{\leq}n)$ of objects of $\ca V$, one obtains a cartesian transformation $\overline{\alpha}:\overline{A}{\rightarrow}E$. The stability of $\ca V$'s coproducts applied to the pullbacks
\[ \PbSq {\overline{A}_n(1,...,1)} {E_n(1,...,1)} {\Gamma{E}(1)} {A1} {} {\overline{\alpha}} {c_n} {\alpha_1} \]
for each $X \in \ca V$ and $n \in \N$ ensures that $\Gamma_E(\overline{\alpha}){\iso}\alpha$ thus verifying essential surjectivity. Let $\alpha:A{\rightarrow}E$ and $\beta:B{\rightarrow}E$ be cartesian, and $\phi:\Gamma{A}{\rightarrow}\Gamma{B}$ such that $\Gamma(\beta)\phi{=}\Gamma{\alpha}$. To finish the proof we must show there is a unique $\phi':A{\rightarrow}B$ such that $\beta\phi'{=}\alpha$ and $\Gamma{\phi'}{=}\phi$. The equation $\Gamma{\phi'}{=}\phi$ implies in particular that $\coprod\limits_n\phi'_{n,1}{=}\phi_1$, and this determines the components $\phi'_{n,1,...,1}$ uniquely because of
%%%%%%%%
\[ \xygraph{!{0;(2.5,0):(0,.5)::}
{A_n(1,...,1)}="a" [r] {B_n(1,...,1)}="b" [r] {E_n(1,...,1)}="e"
[d] {\Gamma{E}(1)}="ge" [l] {\Gamma{B}(1)}="gb" [l] {\Gamma{A}(1)}="ga"
"a" (:@{.>}"b"_-{\phi'_{n,1,...,1}} :"e"_-{\beta_{n,1,...,1}}, :@/^{1.5pc}/"e"^-{\alpha_{n,1,...,1}})
"ga" (:"gb"^-{\phi_1} :"ge"^-{\coprod\beta_n}, :@/_{1.5pc}/"ge"_-{\coprod\alpha_n})
"a":"ga"_{c_n} "b":"gb"^{c_n} "e":"ge"^{c_n}
"a" [d(.3)r(.1)] :@{-} [r(.1)] :@{-} [u(.15)]
"b" [d(.3)r(.1)] :@{-} [r(.1)] :@{-} [u(.15)]}  \]
%%%%%%%%
and these components determine $\phi'$ uniquely because of
%%%%%%%%
\[ \xygraph{!{0;(2.5,0):(0,.5)::}
{\opAA\limits_iX_i}="a" [r] {\opB\limits_iX_i}="b" [r] {\opE\limits_iX_i}="e"
[d] {\opAA\limits_i1}="ge" [l] {\opB\limits_i1}="gb" [l] {\opE\limits_i1}="ga"
"a" (:@{.>}"b"_-{\phi'_{X_i}} :"e"_-{\beta_{X_i}}, :@/^{1.5pc}/"e"^-{\alpha_{X_i}})
"ga" (:"gb"^-{\phi'_{n,1,...,1}} :"ge"^-{\beta_{n,1,...,1}}, :@/_{1.5pc}/"ge"_-{\alpha_{n,1,...,1}})
"a":"ga"_{\opAA\limits_it_{X_i}} "b":"gb"_{\opB\limits_it_{X_i}} "e":"ge"^{\opE\limits_itX_i}
"a" [d(.3)r(.1)] :@{-} [r(.1)] :@{-} [u(.15)]
"b" [d(.3)r(.1)] :@{-} [r(.1)] :@{-} [u(.15)]}  \]
%%%%%%%%
and the equation $\beta\phi'{=}\alpha$. To see that $\Gamma{\phi'}{=}\phi$, that is $\coprod\limits_n\phi'_{n,X,...,X}{=}\phi_X$ for all $X \in \ca V$, one deduces that the inner square in
%%%%%%%%
\[ \xygraph{!{0;(2.5,0):(0,.4)::}
{A_n(X,...,X)}="a" [r] {B_n(X,...,X)}="b" [d] {\Gamma{B}(X)}="gb" [l] {\Gamma{A}(X)}="ga"
[dl] {\Gamma{A}(1)}="ga1" [u(3)] {A_n(1,...,1)}="a1" [r(3)] {B_n(1,...,1)}="b1" [d(3)] {\Gamma{B}(1)}="gb1"
"a" (:"b"^-{\phi'_{n,X,...,X}} :"gb"^{c_n}, :"ga"_{c_n} :"gb"_-{\phi_X})
"a1" (:"b1"^-{\phi'_{n,1,...,1}} :"gb1"^{c_n}, :"ga1"_{c_n} :"gb1"_-{\phi_1})
"a":"a1"|{A_n(t_X,...,t_X)} "b":"b1"|{B_n(t_X,...,t_X)} "ga":"ga1"|{\Gamma{A}(t_X)} "gb":"gb1"|{\Gamma{B}(t_X)}}  \]
%%%%%%%%
is a pullback since the outer square and all other regions in this diagram are pullbacks, and so the result follows by lextensivity.
\end{proof}
%

%%%%%%%%%%%%%%%%%%%%%%%%%%%%%%%%%%%%%%%%%%%%%%%%%%
%%%%%%%%%%%%%%%%%%%%%%%%%%%%%%%%%%%%%%%%%%%%%%%%%%
\section{The strict $\omega$-category monad}\label{sec:wcat-monad}

The setting of the previous section involved a coproduct preserving p.r.a monad $T$, and after this section we shall be concerned with the case where $T=\ca T$ the strict $\omega$-category monad on $\PSh {\G}$ the category of globular sets, and its finite dimensional analogues the strict $n$-category monads. We give a precise and purely inductive combinatorial description of $\ca T$ in section(\ref{ssec:wcat-monad}), using some further theory of p.r.a monads on presheaf categories which we develop in section(\ref{ssec:spra}), to facilitate our description of the details.

%%%%%%%%%%%%%%%%%%%%%%%%%%%%%%%%%%%%%%%%%%%%%%%%%%
\subsection{Specifying p.r.a monads on presheaf categories}\label{ssec:spra}

From \cite{Fam2fun} we know that to specify a p.r.a $T:\PSh {\B}{\rightarrow}\PSh {\C}$ one can begin with $P \in \PSh {\C}$ and a functor $E_T:\el(P){\rightarrow}\PSh {\B}$. Here we will usually not distinguish notationally between $p \in PC$ and $E_T(p,C)$. Given $k:D{\rightarrow}C$ in $\C$ we shall denote by $pk$ the element $Pk(p)$ and by $\overline{k}:pk{\rightarrow}p$ the map
\[ E_T(k:(pk,D){\rightarrow}(p,C)). \]
Given this data one can then define an element of $TX(C)$ to be a pair $(p,h)$ where $p \in PC$ and $h:p{\rightarrow}X$ in $\PSh {\B}$. For a map $k:D{\rightarrow}C$ one defines $TX(k)(p,h)=(pk,h\overline{k})$, and one identifies $P{=}T1$. If the $E_T(p,C)$ are all connected, then $T$ preserves coproducts.

With $T$ so specified it is not hard to characterise generic morphisms. To give a map $f:A{\rightarrow}TX$ is to give for $a \in AC$ an element $p_a \in PC$ together with a map $f_a:p_a{\rightarrow}X$ in $\PSh {\B}$, and this data should be natural in $C$. The assignment $(C,a) \mapsto p_a$ is the object map of a functor $\overline{f}:\el(A){\rightarrow}\PSh {\B}$ and the $f_a$ are the components of a cocone with vertex $X$. Factoring this cocone through its colimit $Z$ gives a factorisation
\[ \xymatrix{A \ar[r]^-{g} & {TZ} \ar[r]^-{Th} & {TX}} \]
where the $g_a$ are the components of the universal cocone. One can easily verify directly that such a $g$ is generic, and since generic factorisations are unique up to isomorphism, one obtains
\begin{lemma}\label{lem:gen-col}
For $T:\PSh {\B}{\rightarrow}\PSh {\C}$ specified as above, $f:A{\rightarrow}TX$ is generic iff its associated cocone exhibits $X$ as a colimit.
\end{lemma} 
\begin{examples}\label{ex:gen}
\begin{enumerate}
\item  If in particular $A$ is a representable $C$, then $f:A{\rightarrow}TX$ amounts to a pair $(p,h:p{\rightarrow}X)$. The associated cocone consists of the one map $p$ and so $f$ is generic in this case iff $p$ is an isomorphism.\label{genex1}
\item  In the case $T=1_{\PSh {\C}}$, $f:A{\rightarrow}X$ is generic iff it is an isomorphism.\label{genex2}
\item  Given $T:\PSh {\C}{\rightarrow}\PSh {\C}$ specified as above, a morphism $f:C{\rightarrow}T^2X$ amounts to a pair $(p,h:p{\rightarrow}TX)$. This morphism is $T^2$-generic iff $h$ is $T$-generic because to give a commuting diagram as depicted on the left
\[ \TwoDiagRel {\xymatrix{{C} \ar[r]^-{\alpha} \ar[d]_{f} & {T^2Y} \ar[d]^{T^2\gamma} \\ {T^2X} \ar[r]_-{T^2\beta} & {T^2Z}}} {} {\xymatrix{{p} \ar[r]^-{\alpha'} \ar[d]_{h} & {TY} \ar[d]^{T\gamma} \\ {TX} \ar[r]_-{T^2\beta} & {TZ}}} \]
is the same as giving a commuting diagram as depicted on the right in the previous display, and so the assertion follows by definition of ``generic''.\label{genex3}
\end{enumerate}
\end{examples}
Suppose now that such a $T:\PSh {\C}{\rightarrow}\PSh {\C}$ comes with a cartesian transformation $\eta:1{\rightarrow}T$. The component $\eta_1$ picks out elements $u_C \in PC$ and for all $X \in \PSh {\C}$ the naturality of $\eta$ with respect to the map $X{\rightarrow}1$ shows that the components of $\eta$ have the explicit form
\[ x \in XC \mapsto (u_C,x':u_C{\rightarrow}X). \]
Observing
\[ \TwoDiagRel {\xymatrix{{u_CC} \ar[r]^-{\eta} \ar[d]_{x'} & {Tu_C(C)} \ar[d]^{Tx'} \\ {XC} \ar[r]_-{\eta} & {TX(C)}}} {} {\xymatrix{{\iota} \ar@{|.>}[r] \ar@{|.>}[d] & {(u_C,1_{u_C})} \ar@{|->}[d] \\ {x} \ar@{|->}[r] & {(u_C,x')}}} \]
we have a unique element of $\iota \in u_CC$ which is sent by $\eta$ to $1_{u_C}$. It is a general fact \cite{WebGen} that components of cartesian transformations reflect generic morphisms, and so by examples(\ref{ex:gen})(\ref{genex1}) and (\ref{genex2}) the morphism $C{\rightarrow}u_C$ corresponding to $\iota$ is an isomorphism. One may assume that this isomorphism is an identity by redefining the yoneda embedding if necessary to agree with $C \mapsto u_C$ and similarly on arrows, so we shall write $C=u_C$. Then the components of $\eta$ may be written as
\[ x \mapsto (C,x:C{\rightarrow}X) \]
where the $x$ on the right hand side corresponds to the $x$ on the left hand side by the yoneda lemma.
\begin{definition}\label{def:pra-spec}
Let $T$ be a p.r.a endofunctor of $\PSh {\C}$ and $\eta:1{\rightarrow}T$ be a cartesian transformation. A pair $(P,E_T)$ giving the explicit description of $(T,\eta)$ as above is called a \emph{specification} of $(T,\eta)$.
\end{definition}
\noindent By the discussion preceeding definition(\ref{def:pra-spec}) every such $(T,\eta)$ has a specification. Let us denote the assignments of an arbitary natural transformation $\mu:T^2{\rightarrow}T$ by
\[ (p \in PC, f:p{\rightarrow}TX) \mapsto (q_f \in PC, h_f:q_f{\rightarrow}X). \]
Naturality of $\mu$ in $C$ says that for $k:D{\rightarrow}C$, $q_{f\overline{k}}=q_fk$ and $h_{f\overline{k}}=h_f\overline{k}$. Naturality of $\mu$ in $X$ says that for $h:X{\rightarrow}Y$, $q_{T(h)f}=q_f$ and $h_{T(h)f}=hh_f$. Suppose that $\mu$ is cartesian. Observing
\[ \TwoDiagRel {\xymatrix{{T^2q_f} \ar[r]^-{\mu} \ar[d]_{T^2h_f} & {Tq_f} \ar[d]^{Th_f} \\ {T^2X} \ar[r]_-{\mu} & {TX}}} {} {\xymatrix{{g_f} \ar@{|.>}[r] \ar@{|.>}[d] & {(q_f,1_{q_f})} \ar@{|->}[d] \\ {(p,f:p{\rightarrow}X)} \ar@{|->}[r] & {(q_f,h_f)}}} \]
one finds that $\forall p \in PC$ and $f:p{\rightarrow}X$, $\exists{!} g_f:p{\rightarrow}Tq_f$ such that $f{=}T(h_f)g_f$ and $h_{g_f}={\id}$. By example(\ref{ex:gen})(\ref{genex3}) and the fact that cartesian transformations reflect generics, such $g_f$'s are automatically generic. Conversely given such $g_f$'s one can readily verify that the naturality squares of $\mu$ corresponding to maps $X{\rightarrow}1$ are pullbacks and so verify that $\mu$ is cartesian. We record these observations in
\begin{lemma}\label{lem:spec-cart-mu}
Let $(T,\eta)$ be specified as in definition(\ref{def:pra-spec}). To give a cartesian natural transformation $\mu:T^2{\rightarrow}T$ is to give for each $p \in PC$ and $f:p{\rightarrow}TX$, an element $q_f \in PC$ and a factorisation
\[ \xymatrix{p \ar[r]^-{g_f} & {Tq_f} \ar[r]^-{Th_f} & {TX}} \]
satisfying
\begin{enumerate}
\item  For $k:D{\rightarrow}C$, $q_{f\overline{k}}=q_fk$ and $h_{f\overline{k}}=h_f\overline{k}$.\label{spec1}
\item  For $h:X{\rightarrow}Y$, $q_{T(h)f}=q_f$ and $h_{T(h)f}=hh_f$.\label{spec2}
\item  For all $p \in PC$ and $f:p{\rightarrow}TX$, $g_f$ is unique such that $f{=}T(h_f)g_f$ and $h_{g_f}=\id$.\label{spec3}
\end{enumerate}
and given this data, the $g_f$ are automatically generic morphisms.
\end{lemma}
\noindent Thus a cartesian transformation $\mu:T^2{\rightarrow}T$ amounts to a nice choice of certain generic factorisations for $T$. Given such a characterisation it is straight-forward to unpack what the monad axioms for $(T,\eta,\mu)$ say in terms of these factorisations.
\begin{lemma}\label{lem:spec-pra-monad}
Let $(T,\eta)$ be specified as in definition(\ref{def:pra-spec}). To give $\mu:T^2{\rightarrow}T$ making $(T,\eta,\mu)$ a p.r.a monad is to give factorisations as in lemma(\ref{lem:spec-cart-mu}) which satisfy the following further conditions:
\begin{enumerate}
\item  For all $p \in PC$ and $f:p{\rightarrow}X$, $q_{\eta{f}}=p$ and $h_{\eta{f}}=f$.
\item  For all $p \in PC$ and $f:p{\rightarrow}X$, $q_{(p,f)}=p$ and $h_{(p,f)}=f$ where $(p,f)$ denotes the map $C{\rightarrow}TX$ corresponding to the element $(p,f) \in TX(C)$ by the yoneda lemma.
\item  For all $p \in PC$ and $f:p{\rightarrow}TX$, $q_{h_f}=q_{\mu{f}}$ and $h_{h_f}=h_{\mu{f}}$.
\end{enumerate}
\end{lemma}
\noindent To summarise, given a specification of a p.r.a $T:\PSh {\C}{\rightarrow}\PSh {\C}$, one has for each $C \in \C$ and $f:C{\rightarrow}TX$, $p \in PC$ and a generic factorisation 
\[ \xymatrix{C \ar[r]^-{g} & {Tp} \ar[r]^-{Th} & {TX}} \]
of $f$. The data of a p.r.a monad $(T,\eta,\mu)$ enables us to regard $C \in PC$ and gives us for each $p \in PC$ and $f:p{\rightarrow}TX$, a choice of $q_f \in PC$ and generic factorisation
\[ \xymatrix{p \ar[r]^-{g_f} & {Tq_f} \ar[r]^-{Th_f} & {TX}} \]
of $f$, and these choices satisfy certain axioms.

In the case of the strict $\omega$-category monad below some further simplifications are possible enabling one to dispense with need to verify the additional conditions of lemma(\ref{lem:spec-pra-monad}) when describing it. The reason as we shall see, is that this case conforms to the following definition.
\begin{definition}\label{def:tight}
A p.r.a $T:\PSh {\B}{\rightarrow}\PSh {\C}$ specified by $E_T:\el(P){\rightarrow}\PSh {\B}$ is \emph{tight} when for all $p$ and $q \in PC$ and $\iota:p \iso q$ in $\PSh {\B}$, one has $p=q$ in $PC$ and $\iota=\id$.
\end{definition}
\noindent Clearly tightness is a property of $T$, that is, is independent of the specification.
\begin{examples}\label{ex:tight}
\begin{enumerate}
\item  Let $T$ be the free monoid endofunctor of $\Set$. Then $E_T$ is a functor $\N{\rightarrow}\Set$ sending $n \in \N$ to a set with $n$ elements. There are of course many non-trivial automorphisms of a finite set, and so $T$ is not tight.
\item  Let $T$ be the free category endofunctor on $\Graph$ which we regard as presheaves on
$\xymatrix@1{0 \ar@<1ex>[r] \ar@<-1ex>[r] & 1}$. Then $P0=\{0\}$ and $P1=\{[n]:n \in \N\}$ and the graph $[n]$ has object set $\{i:0{\leq}i{\leq}n\}$ and a unique edge $(i-1){\rightarrow}i$ for each $1{\leq}i{\leq}n$. With these details at hand one readily verifies that this $T$ is tight.
\item  The free symmetric multicategory endofunctor on the category of multigraphs as described in example(2.14) of \cite{WebGen} is not tight. In this case one actually has distinct $p$ and $q$ in $PC$ sent by $E_T$ to isomorphic multigraphs.
\end{enumerate}
\end{examples}
\begin{lemma}\label{lem:tight->unique}
If $T:\PSh {\B}{\rightarrow}\PSh {\C}$ is a tight p.r.a then for all $A:\PSh {\B}{\rightarrow}\PSh {\C}$ there exists at most one cartesian transformation $A{\rightarrow}T$.
\end{lemma}
\begin{proof}
Let $\alpha$ and $\beta:A{\rightarrow}T$ be cartesian transformations and $a \in AX(C)$. For a given specification $P$ one has $p_{\alpha}$ and a generic factorisation
\[ \xymatrix{C \ar[r]^-{a} \ar[dr]_{g_{\alpha}} & {AX} \ar[r]^-{\alpha_X} & {TX} \\ & {Tp_{\alpha}} \ar[ur]_{Th_{\alpha}}} \]
and using the cartesian naturality square for $\alpha$ corresponding to $h_{\alpha}$, one has $g'_{\alpha}:C{\rightarrow}Ap_{\alpha}$ unique such that $\alpha{g'_{\alpha}}=g_{\alpha}$ and $a=A(h_{\alpha})g'_{\alpha}$. Since cartesian transformations reflect generics, this last equation is an $A$-generic factorisation of $a$, and similarly one obtains another one: $a=A(h_{\beta})g'_{\beta}$ by using $\beta$ instead of $\alpha$. Thus there is a unique isomorphism $\delta:p_{\alpha}{\rightarrow}p_{\beta}$ so that $A(\delta)g'_{\alpha}=g'_{\beta}$ and $h_{\alpha}\delta=h_{\beta}$. By tightness $\delta$ is an identity and so $\alpha_X{a}=\beta_X{a}$.
\end{proof}
\noindent Thus given a tight p.r.a $T:\PSh {\C}{\rightarrow}\PSh {\C}$, cartesian transformations $\eta:1{\rightarrow}T$ and $\mu:T^2{\rightarrow}T$ are unique if they exist, and when they do the monad axioms for $(T,\eta,\mu)$ are automatic. This gives the following refinement of lemma(\ref{lem:spec-pra-monad}) in the tight case.
\begin{corollary}\label{cor:tight-monad}
Let $(T,\eta)$ be specified as in definition(\ref{def:pra-spec}) and let $T$ be tight. To give $\mu:T^2{\rightarrow}T$ making $(T,\eta,\mu)$ a p.r.a monad is to give factorisations as in lemma(\ref{lem:spec-cart-mu}).
\end{corollary}
\noindent Moreover for a tight p.r.a monad $T$ on $\PSh {\C}$, the multitensor $T^{\times}$ admits the same simplifications.
\begin{lemma}\label{lem:Tcross-tight}
Let $(T,\eta,\mu)$ be a p.r.a monad on $\PSh {\C}$ such that $T$ is tight. Then for all $E:\ca M\PSh {\C}{\rightarrow}\PSh {\C}$, there exists at most one cartesian transformation $\varepsilon:E{\rightarrow}T^{\times}$.
\end{lemma}
\begin{proof}
To give such an $\varepsilon$ is to give for each $n \in \N$ a cartesian transformation $\varepsilon_n:E_n{\rightarrow}T^{\times}_n$, and so it suffices by lemma(\ref{lem:tight->unique}), to show that $T^{\times}_n:{\PSh \C}^n{\rightarrow}\PSh {\C}$ is tight for all $n \in \N$. The functor $E_{T^{\times}_n}$ has object map $((p_1,...,p_n),C) \mapsto (p_1,...,p_n)$. For $q_1,...,q_n \in T1(C)$, to give an isomorphism $\iota:(p_1,...,p_n)\iso(q_1,...,q_n)$ in $\PSh {\C}^n$, is to give isomorphisms $\iota_i:p_i\iso{q_i}$ for $1{\leq}i{\leq}n$, in which case the $\iota_i$ are identities by the tightness of $T$, and so $T^{\times}_n$ is also tight.
\end{proof}
\noindent Thus for a tight monad $T$ on $\PSh {\C}$, being a $T$-operad is actually a \emph{property} of a monad on $\PSh {\C}$, and similarly for $T$-multitensors. We shall exploit this observation notationally below, for instance, by denoting a $T$-operad $\alpha:A{\rightarrow}T$ as we just have as a monad morphism, or just by referring to the monad $A$, depending on what is most convenient for the given situation.

%%%%%%%%%%%%%%%%%%%%%%%%%%%%%%%%%%%%%%%%%%%%%%%%%%
\subsection{Inductive description of the strict $\omega$-category monad}\label{ssec:wcat-monad}

A goal of this paper to clarify the inductive nature of the operadic approach to higher category theory of \cite{Bat98}. The starting point of that approach is a precise description of the monad $(\ca T,\eta,\mu)$ on the category $\PSh {\G}$ of globular sets whose algebras are strict $\omega$-categories. Thus in this section we recall this monad, but describe it a little differently to the way it has been described in the past. We shall give here a purely inductive description of this fundamental object, and we shall use the results of the previous section to expedite our account of the details. That the algebras for the monad described in this section really are strict $\omega$-categories defined in the usual way by successive enrichments, is presented in section(\ref{sec:fin-dim}) as a pleasant application of our general theory.

The category $\G$ has as objects natural numbers and for $n<m$ maps
\[ \xymatrix{n \ar@<1ex>[r]^-{\sigma}  \ar@<-1ex>[r]_-{\tau} & m} \]
and these satisfy $\sigma\tau{=}\tau\tau$ and $\tau\sigma{=}\sigma\sigma$.
Thus an object of the category $\PSh {\G}$ of globular sets is a diagram
\[ \xymatrix{{X_0} & {X_1} \ar@<1ex>[l]^{t}  \ar@<-1ex>[l]_{s} & {X_2} \ar@<1ex>[l]^{t}  \ar@<-1ex>[l]_{s}
& {X_3} \ar@<1ex>[l]^{t}  \ar@<-1ex>[l]_{s} & {...} \ar@<1ex>[l]^{t}  \ar@<-1ex>[l]_{s}} \]
of sets and functions such that $ss=st$ and $ts=tt$.
The elements of $X_n$ are called n-cells, and for an $(n+1)$-cell $x$,
the n-cells $sx$ and $tx$ are called the source and target of $x$ respectively.
In fact for each $k{\leq}n$, we can define source and target k-cells of $x$
and we denote these by $s_kx$ and $t_kx$, only dropping the indexing
when there is little risk of confusion. 
Given a pair $(a,b)$ of $n$-cells of $X$, one can define the globular set
$X(a,b)$. A $k$-cell of $X(a,b)$ is an $(n+k)$-cell $x$ of $X$
such that $s_kx=a$ and $t_kx=b$. Sources and targets for $X(a,b)$
are inherited from $X$. In particular the globular sets $X(a,b)$ where $a$ and $b$ are 0-cells are called \emph{the homs of $X$}. A morphism $f:X{\rightarrow}Z$ of globular sets
induces maps $X(a,b){\rightarrow}Z(f_0a,f_0b)$ on the homs.
Conversely, to give $f$ it suffices to specify a function $f_0:X_0{\rightarrow}Z_0$
and for all $a, b \in X_0$, morphisms $X(a,b){\rightarrow}Z(f_0a,f_0b)$
of globular sets.

A finite sequence $(X_1,...,X_n)$ of globular sets may be regarded as a globular set, whose set of $0$-cells is $\{i \in \N : 0 \leq i \leq n\}$ and whose only non-empty homs are given by $(X_1,...,X_n)(i-1,i) = X_i$ for $1{\leq}i{\leq}n$. This construction is the object map of a functor $\PSh {\G}^n{\rightarrow}\PSh {\G}$.

We now begin our description of the endofunctor $\ca T$ in the spirit of section(\ref{ssec:spra}). The role of $P$ is played by the globular set $\Tr$ of trees. The set $\Tr_0$ contains one element denoted as $0$ and its associated globular set contains one $0$-cell, also called $0$, and nothing else. By induction an element of $\Tr_{n+1}$ is a finite sequence $(p_1,...,p_k)$ of elements of $\Tr_n$ and its associated globular set is just the sequence of globular sets $(p_1,...,p_k)$ regarded as a globular set as in the previous paragraph. So far we have defined the elements of $\Tr_n$ for all $n$ and the object map of $E_{\ca T}:\el(\Tr){\rightarrow}\PSh {\G}$. We denote by $\sigma:0{\rightarrow}p$ the map which selects the object $0 \in p$, and by $\tau:0{\rightarrow}p$ the map which selects the maximum vertex of $p$ (using $\leq$ inherited from $\N$).

The source and target maps $s,t:\Tr_{n+1}{\rightarrow}\Tr_n$ coincide and are denoted as $\partial$. For each $n$ we must define this map and give maps $\sigma:\partial{p}{\rightarrow}p$ and $\tau:\partial{p}{\rightarrow}p$ which satisfy the equations $\sigma\sigma=\tau\sigma$ and $\tau\tau=\sigma\tau$ in
\begin{equation}\label{eq:csct} \xymatrix{{\partial^2p} \ar@<1ex>[r]^-{\sigma} \ar@<-1ex>[r]_-{\tau} & {\partial{p}} \ar@<1ex>[r]^-{\sigma} \ar@<-1ex>[r]_-{\tau} & {p}} \end{equation}
for all $p \in \Tr_{n+2}$, in order to complete the description of $\Tr$ and the functor $E_{\ca T}$, and thus the definition of $\ca T$. The maps $\partial$, $\sigma$ and $\tau$ are given by induction as follows. For the initial step $\partial$ is uniquely determined since $\Tr_0$ is singleton and $\sigma$ and $\tau$ are as described in the previous paragraph. For the inductive step let $p=(p_1,...,p_k) \in \Tr_{n+2}$. Then $\partial{p}=(\partial{p_1},...,\partial{p_k})$ and the maps $\sigma,\tau:\partial{p}{\rightarrow}p$ are the identities on $0$-cells, and the non-empty hom maps are given by $\sigma,\tau:\partial{p_i}{\rightarrow}p_i$ respectively for $1{\leq}i{\leq}k$. The verification of $\sigma\sigma=\tau\sigma$ and $\tau\tau=\sigma\tau$ as in (\ref{eq:csct}) is given by induction as follows. The initial step when $n=0$ is clear since $\partial^2p=0$, and the $0$-cell maps of $\sigma,\tau:\partial{p}{\rightarrow}p$ are both the identity. For the inductive step let $p !
 \in \Tr_{n+3}$, then all the maps in (\ref{eq:csct}) are identities on $0$-cells, and on the homs the desired equations follow by induction.

By section(\ref{ssec:spra}) we have completed the description of a p.r.a $\ca T:\PSh{\G}{\rightarrow}\PSh {\G}$ and we will now see that it is tight. Once again we argue by induction on $n$. In the case $n=0$ the result follows because $\Tr_0=\{0\}$ and the only automorphism of $0 \in \PSh {\G}$ is the identity. For the inductive step let $p,q \in \Tr_{n+1}$ and suppose that one has $\iota:p \iso q$ in $\PSh {\G}$. Since the only non-empty homs for $p$ and $q$ are between consecutive elements of their vertex sets, any $f:p{\rightarrow}q$ in $\PSh {\G}$ is order preserving in dimension $0$. Thus the $0$-cell map of $\iota$ is an order preserving bijection, and so must be the identity. The hom maps of $\iota$ must also be identities by induction. Since the globular sets associated to $p \in \Tr_n$ are also connected we have the following result.
\begin{proposition}\label{prop:wcat-tightpra}
$\ca T:\PSh {\G}{\rightarrow}\PSh {\G}$ defined as follows is p.r.a, tight and coproduct preserving:
\begin{itemize}
\item  an $n$-cell of $\ca TX$ is a pair $(p,f:p{\rightarrow}X)$ where $p \in \Tr_n$.
\item  for $n{\geq}1$, $s(p,f)=(\partial{p},f\sigma)$ and $t(p,f)=(\partial{p},f\tau)$.
\item  for $h:X{\rightarrow}Y$, $\ca T(h)(p,f)=(p,hf)$.
\end{itemize}
\end{proposition}
We will now specify the cartesian unit $\eta:1{\rightarrow}\ca T$, and from section(\ref{ssec:spra}) we know that this amounts to factoring the yoneda embedding through $E_{\ca T}$. We already have $0 \in \Tr_0$, and by induction we define $n+1=(n) \in \PSh {\G}$. Notice that the set of $k$-cells of $n$ is $\{0,1\}$ when $k<n$ and $\{0\}$ when $k=n$. Moreover by an easy inductive proof the reader may verify that the $k$-cell maps of $\sigma:n{\rightarrow}n+1$ and $\tau:n{\rightarrow}n+1$ are the identities for $k<n$, and pick out $0$ and $1$ respectively when $k=n$. One has functions $\ev_0:\PSh {\G}(n,X){\rightarrow}X_n$ given by $f \mapsto f_n(0)$ clearly natural in $X \in \PSh {\G}$. By another easy induction one may verify that these functions are bijective, and natural in $n$ in the sense that $sf_{n+1}(0)=(f\sigma)_n(0)$ and $tf_{n+1}(0)=(f\tau)_n(0)$. Henceforth we regard the identification of $n$ as a globular set in this way as \emph{the} yoneda embedding, and the c!
 
 omponents of $\eta$ are given by $x \in X_n \mapsto x:n{\rightarrow}X$.

Before specifying the multiplication $\mu:\ca T^2{\rightarrow}\ca T$ some preliminary remarks are in order. For $0$-cells $a$ and $b$ of $X$, an $n$-cell of the hom $\ca TX(a,b)$ consists by definition, of $p=(p_1,...,p_k) \in \Tr_{n+1}$ together with $f:p{\rightarrow}X$ such that $f\sigma=a$ and $f\tau=b$. In other words one has a sequence $(x_0,...,x_k)$ of $0$-cells of $X$ such that $x_0=a$ and $x_k=b$, together with maps $f_i:p_i{\rightarrow}X(x_{i-1},x_i)$ for $1{\leq}i{\leq}k$. Another way to say all this is that for a given sequence $(x_0,...,x_k)$ of $0$-cells of $X$ such that $x_0=a$ and $x_k=b$, one has an inclusion
\[ \begin{array}{c} {c_{x_i} : \prod\limits_{1{\leq}i{\leq}k} \ca T(X(x_{i-1},x_i)) \rightarrow \ca TX(a,b)} \end{array} \]
in $\PSh {\G}$, and the following result.
\begin{lemma}\label{lem:homs-of-TX}
The maps $c_{x_i}$, for all sequences $(x_0,...,x_k)$ of $0$-cells of $X$ such that $x_0=a$ and $x_k=b$, form a coproduct cocone.
\end{lemma}
\noindent Let $p=(p_1,...,p_k) \in \Tr_{n+1}$. A map $f:p{\rightarrow}\ca TX$ amounts to $0$-cells $fi$ of $X$ for $0{\leq}i{\leq}k$, together with hom maps $f_i:p_i{\rightarrow}\ca TX(f(i-1),fi)$ for $1{\leq}i{\leq}k$. Since the $p_i$ are connected, the $f_i$ amount to $0$-cells $(x_{i0},...,x_{im_i})$ of $X$ such that $x_{i0}=f(i-1)$ and $x_{im_i}=fi$, together with maps $f_{ij}:p_i{\rightarrow}\ca T(X(x_{(ij)-1},x_{ij}))$ for $1{\leq}i{\leq}k$ and $1{\leq}j{\leq}m_i$ where
\[ (i,j)-1 = \left\{\begin{array}{lll} {(i,j-1)} && {\textnormal{when $j>0$.}} \\ {(i-1,m_{i-1})} && {\textnormal{when $j=0$ and $i>0$.}} \\ {0} && {\textnormal{when $i=j=0$.}} \end{array}\right. \]
In other words for $p=(p_1,...,p_k) \in \Tr_{n+1}$, to give $f:p{\rightarrow}\ca TX$ is to give objects $x_0$ and $x_{ij}$ of $X$ together with maps $f_{ij}:p_i{\rightarrow}\ca T(X(x_{(ij)-1},x_{ij}))$ for $1{\leq}i{\leq}k$ and $1{\leq}j{\leq}m_i$. We shall call $x_0$ and the $x_{ij}$ \emph{the $0$-cells of $f$}, and the $f_{ij}$ the \emph{hom map components} of $f$. Observe that for $h:X{\rightarrow}Y$, the $0$-cells of $\ca T(h)f$ are given by $hx_0$ and $hx_{ij}$, and the hom map components by $hf_{ij}$ where $1{\leq}i{\leq}k$ and $1{\leq}j{\leq}m_i$.

Now we specify the multiplication $\mu:\ca T^2{\rightarrow}\ca T$ following lemma(\ref{lem:spec-cart-mu}). For $p \in \Tr_n$ and $f:p{\rightarrow}\ca TX$ the factorisation of $f$ that we must provide will be given by induction on $n$. When $n=0$, $p=0$ and a map $f:0{\rightarrow}TX$ picks out a $0$-cell $(0,x:0{\rightarrow}X)$ of $\ca TX$. Define $q_f=0$, $h_f=x$ and $g_f:0{\rightarrow}\ca T0$ to pick out $(0,1_0)$. For the inductive step let $p=(p_1,...,p_k) \in \Tr_{n+1}$ and $f:p{\rightarrow}\ca TX$. Then define
\[ q_f = (q_{f_{ij}} : 1{\leq}i{\leq}k, 1{\leq}j{\leq}m_i) \]
where the $f_{ij}$ are the hom map components of $f$ as defined in the previous paragraph. Define $h_f$ to have $0$-cell mapping given by $0{\mapsto}x_0$ and $(i,j){\mapsto}x_{ij}$, and hom maps by $h_{f_{ij}}$. Define $g_f$ to have underlying $0$-cells given by $0$ and $(i,j)$, and hom map components by $g_{f_{ij}}$. By definition we have $f=\ca T(h_f)g_f$.
\begin{proposition}\label{prop:wcat-pra-monad}
$(\ca T,\eta,\mu)$ with $\ca T$ as specified in proposition(\ref{prop:wcat-tightpra}), and $\eta$ and $\mu$ given by
\[ \begin{array}{lccr} {x \in X_n \mapsto (n,x:n{\rightarrow}X)} &&& {(p \in \Tr_n, f:p{\rightarrow}\ca TX) \mapsto (q_f,h_f)} \end{array} \]
is a p.r.a monad.
\end{proposition}
\begin{proof}
By corollary(\ref{cor:tight-monad}) it suffices to verify conditions (\ref{spec1})-(\ref{spec3}) of lemma(\ref{lem:spec-cart-mu}). Condition(\ref{spec1}) says that for $p=(p_1,...,p_k) \in \Tr_{n+1}$ and $f:p{\rightarrow}\ca TX$: $q_{f\sigma}=q_{f\tau}=\partial{q_{f}}$, $h_{f\sigma}=h_f\sigma$ and $h_{f\tau}=h_f\tau$. Let us write $x_0$ and $x_{ij}$ for the $0$-cells of $f$ and $f_{ij}$ for the hom map components where $1{\leq}i{\leq}k$ and $1{\leq}j{\leq}m_i$. In the case $n=0$, we must have $0=q_{f\sigma}=q_{f\tau}=\partial{q_{f}}$ since $0$ is the only element of $\Tr_0$. Clearly $f\sigma$ picks out $x_0$ and $f\tau$ picks out $x_{km_k}$, and so $h_{f\sigma}:0{\rightarrow}X$ picks out $x_0$ and $h_{f\tau}:0{\rightarrow}X$ picks out $x_{km_k}$ by the initial step of the description of the factorisations. By the definition of the object map of $h_f$, $h_f\sigma$ and $h_f\tau$ also pick out the $0$-cells $x_0$ and $x_{km_k}$ respectively, thus verifying the $n=0$ case of con!
 
 dition(\ref{spec1}). For the inductive step let $p=(p_1,...,p_k) \in \Tr_{n+2}$ and $f:p{\rightarrow}\ca TX$. First note that $\sigma,\tau:\partial{p}{\rightarrow}p$ are identities on $0$-cells and so $f$, $f\sigma$ and $f\tau$ have the same $0$-cells which we are denoting by $x_0$ and $x_{ij}$. Moreover by the definition of hom map components, one has $(f\sigma)_{ij}=f_{ij}\sigma$ and $(f\tau)_{ij}=f_{ij}\tau$. Thus by induction
\[ q_{f\sigma} = (q_{f_{ij}\sigma} : 1{\leq}i{\leq}k, 1{\leq}j{\leq}m_i) = (\partial{q_{f_{ij}}} : 1{\leq}i{\leq}k, 1{\leq}j{\leq}m_i) = \partial{q_f} \]
and similarly $q_{f\tau}=\partial{q_f}$. Since $\sigma,\tau:\partial{q_f}{\rightarrow}q_f$ are identities on $0$-cells the equations $h_{f\sigma}=h_f\sigma$ and $h_{f\tau}=h_f\tau$ are true on $0$-cells, and on homs these equations follow by induction.

Condition(\ref{spec2}) says that for $p \in \Tr_n$, $f:p{\rightarrow}\ca TX$ and $h:X{\rightarrow}Y$, $q_{\ca T(h)f}=q_f$ and $h_{\ca T(h)f}=hh_f$. When $n=0$ these equations are immediate. For the inductive step let $p=(p_1,...,p_k) \in \Tr_{n+1}$, $f:p{\rightarrow}\ca TX$ and $h:X{\rightarrow}Y$. The objects of $q_{\ca T(h)f}$ and $q_f$ coincide by definition, and the homs do by induction. The object maps of $h_{\ca T(h)f}$ and $hh_f$ coincide by definition and their homs maps coincide by induction.

Condition(\ref{spec3}) says that for $p \in \Tr_n$ and $f:p{\rightarrow}\ca TX$, $g_f$ is unique such that $f=\ca T(h_f)g_f$ and $h_{g_f}=\id$. For $n=0$ this is clear by inspection. For the inductive step let $p=(p_1,...,p_k) \in \Tr_{n+1}$ and $f:p{\rightarrow}\ca TX$. By inspection the $0$-cell map of $h_{g_f}$ is the identity, and by induction its hom maps are also identities. As for uniqueness, the object map of $g_f$ is determined uniquely by $k$ and $m_i \in \N$ for $1{\leq}i{\leq}k$, and the uniqueness of the hom maps follows by induction.
\end{proof}
%

%%%%%%%%%%%%%%%%%%%%%%%%%%%%%%%%%%%%%%%%%%%%%%%%%%
%%%%%%%%%%%%%%%%%%%%%%%%%%%%%%%%%%%%%%%%%%%%%%%%%%
\section{Normalised $\ca T$-operads and $\ca T$-multitensors}\label{ssec:corr}
In this section we relate $\ca T$-operads to $\ca T$-multitensors and so express $\ca T$-operad algebras as enriched categories. Under a mild condition on an operad $\alpha:A{\rightarrow}\ca T$, that it be \emph{normalised} in the sense to be defined shortly, one can construct a multitensor $\overline{A}$ on $\PSh {\G}$ such that $\overline{A}$-categories are $A$-algebras. Moreover $\overline{A}$ is in fact a $\ca T$-multitensor, and the construction $\overline{(\,\,)}$ is part of an equivalence of categories between $\Mult {\ca T}$ and the full subcategory of $\Op {\ca T}$ consisting of the normalised $\ca T$-operads.
\begin{definition}
An endofunctor $A$ of $\PSh {\G}$ is \emph{normalised} when for all $X \in \PSh {\G}$, $\{AX\}_0{\iso}X_0$. A monad $(A,\eta,\mu)$ is normalised when $A$ is normalised as an endofunctor, a cartesian transformation $\alpha:A{\rightarrow}\ca T$ is called a \emph{normalised collection} when $A$ is normalised, and a $\ca T$-operad $\alpha:A{\rightarrow}\ca T$ is normalised when $A$ is normalised as a monad or endofunctor. We shall denote by $\NColl {\ca T}$ the full subcategory of $\PraEnd(\PSh {\G})/\ca T$ consisting of the normalised collections, and by $\NOp {\ca T}$ the full subcategory of $\Op {\ca T}$ consisting of the normalised operads.
\end{definition}
\noindent A $0$-cell of $\ca TX$ is a pair $(p \in \Tr_0, x:p{\rightarrow}X)$, but then $p=0$ and by the yoneda lemma we can regard $x$ as an element of $X_0$. Thus $\ca T$ is normalised. The category $\NColl {\ca T}$ inherits a strict monoidal structure from $\PraEnd(\PSh {\G})/\ca T$, and the category of monoids therein is exactly $\NOp {\ca T}$. We shall allow a very convenient abuse of notation and language: for normalised $A$ write $\{AX\}_0{=}X_0$ rather than acknowledging the bijection, and speak of $X$ and $AX$ as having the \emph{same} $0$-cells. This abuse is justified because for any normalised $A$, one can obviously redefine $A$ to $A'$ which is normalised in this strict sense, and the assignment $A \mapsto A'$ is part of an equivalence of categories between normalised endofunctors and ``strictly normalised endofunctors'', regarded as full subcategories of $\End(\PSh {\G})$.

We begin by recalling and setting up some notation. Recall how a finite sequence $(X_1,...,X_k)$ of globular sets may be regarded as a globular set: the set of $0$-cells is
\[ [k]_0 = \{0,...,k\}, \]
$(X_1,...,X_k)(i-1,i)=X_i$ and all the other homs are empty. Since we shall use these sequences often thoughout this section it is necessary to be careful with the use of round brackets with globular sets. For instance $X$ and $(X)$ are different, and so for an endofunctor $A$ of $\PSh {\G}$, one cannot identify $AX$ and $A(X)$!! Observe also that the $0$-cell map of a morphism
\[ f : (X_1,...,X_m) \rightarrow (Y_1,...,Y_n) \]
must be distance preserving, that is it sends consecutive elements to consecutive elements, whenever all the $X_i$ are non-empty globular sets. We regard sequences $(x_0,...,x_k)$ of $0$-cells of a globular set $X$ as maps $x:[k]_0{\rightarrow}X$ in $\PSh {\G}$. Given any such $x$ we shall define
\[ x^*X := (X(x_{i-1},x_i) \, : \, 1{\leq}i{\leq}k), \]
and a map $\overline{x}:x^*X{\rightarrow}X$ of globular sets. The maps $\overline{x}$ and $x$ agree on $0$-cells, and $\overline{x}_{i-1,i}=\id$ for $1{\leq}i{\leq}k$ specifies the hom maps of $\overline{x}$.

Fundamental to this section is the description of the homs of $\ca TX$ given in lemma(\ref{lem:homs-of-TX}). We shall now refine this and see that an analogous lemma holds for any normalised collection. For $X$ a globular set and $a$ and $b \in X_0$, we shall now understand the hom $\{\ca TX\}(a,b)$. An $n$-cell of $\{\ca TX\}(a,b)$ is a pair $(p,f)$ where $p \in \Tr_{n+1}$ and $f:p{\rightarrow}X$, such that $f\sigma=a$ and $f\tau=b$. First we consider the case $X=(X_1,...,X_k)$ for globular sets $X_i$. Writing $p=(p_1,...,p_m)$ where the $p_i \in \Tr_n$, notice that $f_0$ must be distance preserving. There will be no such $f$ when $a>b$, and in the case $a \leq b$ an $n$-cell of $\{\ca TX\}(a,b)$ consists of $p_i \in \Tr_n$ where $a{<}i{\leq}b$ together with $f_i:p_i{\rightarrow}X_i$. In particular note that when $a=0$ and $b=k$, $f_0=\id$. We record this in the following result.
\begin{lemma}\label{lem:Thom1}
Let $X=(X_1,...,X_k)$ in $\PSh {\G}$. Then for $0{\leq}a,b{\leq}k$ we have
\[ \{\ca TX\}(a,b) = \left\{\begin{array}{lcr} {\emptyset} && {a > b} \\
{\prod\limits_{a{<}i{\leq}b}\ca TX_i} && {a \leq b} \end{array}\right. \]
In particular $\ca T^{\times}X_i=\{\ca TX\}(0,k)$.
\end{lemma}
\noindent Now take $X$ to be an arbitrary globular set. Writing $x:[m]_0{\rightarrow}X$ for the sequence of $0$-cells of $X$ defined by $f_0$, notice that $f:p{\rightarrow}X$ factors uniquely as
\[ \xymatrix{p \ar[r]^-{f'} & {x^*X} \ar[r]^-{\overline{x}} & X} \]
and so defines $(p,f') \in \ca Tx^*X_n$ which gets sent to $(p,f)$ by $\ca T\overline{x}$. Notice that $f'$ is the identity on $0$-cells, which is to say that $(p,f')$ is an $n$-cell of $\{\ca Tx^*X\}(0,m)$.
Therefore an $n$-cell $\phi$ of $\{\ca TX\}(a,b)$ is determined uniquely by the following data: $m \in \N$, $x:[m]_0{\rightarrow}X$ such that $x0=a$ and $xm=b$, and an $n$-cell $\phi'$ of $\{\ca Tx^*X\}(0,m)$. One recovers $\phi$ from this data by $\{\ca T\overline{x}\}_{0,m}\phi'=\phi$. Notice also that if any of the homs $X(x(i-1),x_i)$ is empty, then since $\ca T$ preserves the initial object, one has that $\{\ca Tx^*X\}(0,m)$ is empty by lemma(\ref{lem:Thom1}). Thus one can also specify an $n$-cell $\phi$ of $\{\ca TX\}(a,b)$ uniquely by giving $(m,\phi',x)$ as above, with the additional condition on $x$ that the homs $X(x(i-1),xi)$ be non-empty for all $1{\leq}i{\leq}k$. This last condition amounts to saying that one can factor $x$ as
\[ \xymatrix{{[m]_0} \ar[r]^-{i} & {[m]} \ar[r] & X} \]
where $i$ is the inclusion of the vertices of $[m]$. We shall call the sequences $x$ satisfying this condition \emph{connected}. We have proved the following refinement of lemma(\ref{lem:homs-of-TX}).
\begin{lemma}\label{lem:Thom2}
Let $X$ be a globular set and $a$ and $b \in X_0$.
\begin{enumerate}
\item  The maps
\[ \{\ca T\overline{x}\}_{0,m} : \{\ca Tx^*X\}(0,m) \rightarrow \{\ca TX\}(a,b) \]
for all $m \in \N$ and all sequences $x:[m]_0{\rightarrow}X$ such that $x0=a$ and $xm=b$, form a coproduct cocone.
\item  The maps $\{\ca T\overline{x}\}_{0,m}$ for all $m \in \N$ and all connected sequences $x:[m]_0{\rightarrow}X$ such that $x0=a$ and $xm=b$, form a coproduct cocone.
\end{enumerate}
\end{lemma}
\noindent Now for a normalised collection $\alpha:A{\rightarrow}\ca T$, the extensivity of $\PSh {\G}$ and the cartesianness of $\alpha$ enables us to lift our understanding of the homs of $\ca TX$ expressed in the previous two lemmas, to an understanding of the homs of $AX$. In order to do this in lemma(\ref{lem:A-homs}) below, we require a basic lemma regarding pullbacks and homs in $\PSh {\G}$.
\begin{lemma}\label{lem:pb-hom}
Given a commutative square (I)
%%%%%%
\[ \TwoDiagRel {\xymatrix{W \ar[d]_{f} \save \POS?="d" \restore \ar[r]^-{h} & X \ar[d]^{g} \save \POS?="c" \restore \\ Y \ar[r]_-{k} & Z \POS "d";"c" **@{}; ?*{\textnormal{I}}}} {}
{\xymatrix{{W(a,b)} \ar[d]_{f_{a,b}} \save \POS?="d" \restore \ar[r]^-{h_{a,b}} & {X(ha,hb)} \ar[d]^{g_{ha,hb}} \save \POS?="c" \restore \\ {Y(a,b)} \ar[r]_-{k_{a,b}} & {Z(ha,hb)} \POS "d";"c" **@{}; ?*{\textnormal{II}}}} \]
%%%%%%
in $\PSh {\G}$ such that $f_0$ and $g_0$ are identities, one has for each $a,b \in W_0$ commuting squares (II) as in the previous display. The square (I) is a pullback iff for all $a,b \in W_0$, the square (II) is a pullback.
\end{lemma}
\begin{proof}
Suppose that (I) is a pullback and $a,b \in W_0$. Let $y \in Y(a,b)_n$ and $x \in X(ha,hb)_n$ such that $ky=gx$. Then there is a unique $w \in W_{n+1}$ such that $fw=y$ and $hw=x$, and since $f_0=\id$ and its components commute with sources and targets, one has $w \in W(a,b)_n$ whence (II) is a pullback. Conversely suppose that (II) is a pullback for all $a,b \in W_0$. In dimension $0$ (I) is a pullback since $f_0$ and $g_0$ are identities. For $n \in \N$ let $y \in Y_{n+1}$ and $x \in X_{n+1}$ such that $ky=gx$. Put $a=s_0y$ and $b=t_0b$ so that $y \in Y(a,b)_n$. Since the components of maps in $\PSh {\G}$ commute with sources and targets we have $x \in X(ha,hb)_n$, and since (II) is a pullback there is a unique $w \in W(a,b)_n$ such that $fw=y$ and $hw=x$. Any $w' \in W_{n+1}$ such that $fw'=y$ and $hw'=x$ is in $W(a,b)_n$ since the components of $f$ commute with sources and targets, and so $w'=w$.
\end{proof}
\begin{lemma}\label{lem:A-homs}
Fix a choice of initial object $\emptyset$ and pullbacks in $\PSh {\G}$, such that the pullback of an identity arrow is an identity. Let $\alpha:A{\rightarrow}\ca T$ be a normalised collection.
\begin{enumerate}
\item  Let $X=(X_1,...,X_k)$ in $\PSh {\G}$. Then for $0{\leq}a,b{\leq}k$ we have
\[ \{AX\}(a,b) = \left\{\begin{array}{lcr} {\emptyset} && {a > b} \\
{\{Ax^*X\}(0,b-a)} && {a \leq b} \end{array}\right. \]
where $x:[b-a]_0{\rightarrow}X$ is given by $xi=a+i$.\label{A-homs1}
\item  The maps
\[ \{A\overline{x}\}_{0,m} : \{Ax^*X\}(0,m) \rightarrow \{AX\}(a,b) \]
for all $m \in \N$ and all sequences $x:[m]_0{\rightarrow}X$ such that $x0=a$ and $xm=b$, form a coproduct cocone.\label{A-homs2}
\item  The maps $\{A\overline{x}\}_{0,m}$ for all $m \in \N$ and all connected sequences $x:[m]_0{\rightarrow}X$ such that $x0=a$ and $xm=b$, form a coproduct cocone.\label{A-homs3}
\end{enumerate}
\end{lemma}
\begin{proof}
In the case $X=(X_1,...,X_k)$ with $a>b$ one has $\{\alpha_X\}_{a,b}:\{AX\}(a,b){\rightarrow}\emptyset$ by lemma(\ref{lem:Thom1}), and since the initial object of $\PSh {\G}$ is strict, one has $\{AX\}(a,b)=\emptyset$. Given any $X$ and $x:[m]_0{\rightarrow}X$ such that $x0=a$ and $xm=b$, we have that
\[ \xymatrix{{\{Ax^*X\}(0,m)} \ar[r]^-{\{A\overline{x}\}_{0,m}} \ar[d]_{\{\alpha_{x^*X}\}_{0,m}} & {\{AX\}(a,b)} \ar[d]^{\{\alpha_X\}_{a,b}} \\ {\{\ca Tx^*X\}(0,m)} \ar[r]_-{\{\ca T\overline{x}\}_{0,m}} & {\{\ca TX\}(a,b)}} \]
is a pullback by lemma(\ref{lem:pb-hom}) and the cartesianness of $\alpha$. In the case $X=(X_1,...,X_k)$ with $a \leq b$ and $x:[b-a]_0{\rightarrow}X$ given by $xi=a+i$, $\{\ca T\overline{x}\}_{0,b-a}$ is the identity by lemma(\ref{lem:Thom1}), thus so is $\{A\overline{x}\}_{0,b-a}$ and we have proved (\ref{A-homs1}). In the general case considering all $m \in \N$ and sequences (resp. connected sequences) $x:[m]_0{\rightarrow}X$ with $x0=a$ and $xm=b$, the $\{\ca T\overline{x}\}_{0,m}$ form a coproduct cocone by lemma(\ref{lem:Thom2}), and thus so do the $\{A\overline{x}\}_{0,m}$ by extensivity, which gives (\ref{A-homs2}) and (\ref{A-homs3}).
\end{proof}
For a normalised collection $A$, $k \in \N$ and $X_i \in \PSh {\G}$ where $1 \leq i \leq k$, define
\[ \opA\limits_iX_i = \{AX\}(0,k) \]
where $X=(X_1,...,X_k)$.{\footnotemark{\footnotetext{There is an analogy between lemma(\ref{lem:A-homs}) and the Lagrangian formulation of quantum mechanics. In this analogy one regards any globular set $X$, to which one would apply a collection, as a \emph{state space} the $0$-cells of which are called \emph{states}. A normalised collection $A$ is then a \emph{type of quantum mechanical process}, with the hom $\{AX\}(a,b)$ playing the role of the \emph{amplitude} that the process starts in state $a$ and finishes in state $b$. The \emph{basic} amplitudes are the $\{AX\}(0,k)$ where $X=(X_1,...,X_k)$. In terms of these analogies, lemma(\ref{lem:A-homs}) expresses the sense in which the general amplitude $\{AX\}(a,b)$ may be regarded as the sum of the basic amplitudes over all the ``paths'' between $a$ and $b$, that is, as a sort of discrete Feynman integral. The formula just given expresses this passage between basic and general amplitudes as a particular strong monoidal func!
 
 tor, which allows us to view normalised operads as multitensors, and algebras of such an operad as categories enriched in the corresponding multitensor.

The reader should be aware that it first became apparent to the authors that lemma(\ref{lem:A-homs}) is fundamental to the proof of theorems(\ref{thm:bar-monoidal}) and (\ref{thm:bar-equiv}), and the above analogy was noticed afterwards.}}}
\begin{theorem}\label{thm:bar-monoidal}
The assignment $A \mapsto \overline{A}$ is the object map of a strong monoidal functor
\[ \overline{(\,\,)} : \NColl {\ca T} \rightarrow \Dist(\PSh {\G}). \]
For a normalised operad $A$, one has an isomorphism $\Alg A \iso \Enrich {\overline{A}}$ commuting with the forgetful functors into $\Set$.
\end{theorem}
\begin{proof}
The above definition is clearly functorial in the $X_i$ so one has $\overline{A}:\ca M\PSh {\G}{\rightarrow}\PSh {\G}$. A morphism of normalised collections $\phi:A{\rightarrow}B$ is a cartesian transformation between $A$ and $B$, and such a $\phi$ then induces a natural transformation $\overline{\phi}:\overline{A}{\rightarrow}\overline{B}$ by the formula $\overline{\phi}_{X_i}=\{\phi_X\}_{0,k}$. The cartesianness of $\phi$ and lemma(\ref{lem:pb-hom}) ensures that $\overline{\phi}$ is cartesian. In particular $\overline{\ca T}=\ca T^{\times}$ by lemma(\ref{lem:Thom1}) and so for a given normalised collection $\alpha:A{\rightarrow}\ca T$, one obtains a cartesian $\overline{\alpha}:\overline{A}{\rightarrow}\ca T^{\times}$. Now by example(\ref{ex:Tcross-dist}) $\ca T^{\times}$ is distributive (ie preserves coproducts in each variable) and so $\overline{A}$ is also because of the cartesianness of $\overline{\alpha}$ and the stability of coproducts in $\PSh {\G}$. The assignment !
 
 $\phi \mapsto \overline{\phi}$ described above is clearly functorial, and so $\overline{(\,\,)}$ is indeed well-defined as a functor into $\Dist(\PSh {\G})$.

Since $X(0,k)$ is empty when $k \neq 1$ and just $X_1$ when $k=1$, we have $\overline{1}=I$ the unit of $\Dist(\PSh {\G})$. Let $A$ and $B$ be normalised collections and $X=(X_1,...,X_m)$. By lemma(\ref{lem:A-homs}) the morphisms
\[ \{A\overline{x}\}_{0,k} : \{Ax^*BX\}(0,k) \rightarrow \{ABX\}(0,m) \]
where $k \in \N$ and $x:[k]_0{\rightarrow}BX$ such that $x0=0$ and $xk=m$, form a coproduct cocone. By the definition of the tensor product in $\Dist(\PSh {\G})$, this induces an isomorphism $\overline{AB} \iso \overline{A}\comp\overline{B}$. We now argue that these isomorphisms satisfy the coherence conditions of a strong monoidal functor. Recall that the tensor product in $\Dist(\PSh {\G})$ is defined using coproducts. A different choices of coproducts give rise to different monoidal structures on $\Dist(\PSh {\G})$, though for two such choices the identity functor on $\Dist(\PSh {\G})$ inherits unique coherence isomorphisms that make it strong monoidal and thus an isomorphism of monoidal categories. Because of this one may easily check that if a given strong monoidal coherence diagram commutes for a particular choice of defining coproducts of the monoidal structure of $\Dist(\PSh {\G})$, then this diagram commutes for any such choice. Thus to verify a given strong monoida!
 
 l coherence diagram, it suffices to see that it commutes for some choice of coproducts. But for any such diagram one can simply choose the coproducts so that all the coherence isomorphisms involved in \emph{just} that diagram are identities. Note that this is \emph{not} the same as specifying $\Dist(\PSh {\G})$'s monoidal structure so as to make $\overline{(\,\,)}$ strict monoidal. This finishes the proof that $\overline{(\,\,)}$ is strong monoidal.

Let $A$ be a normalised operad and $Z$ be a set. To give a globular set $X$ with $X_0=Z$ and $x:AX{\rightarrow}X$ which is the identity on $0$-cells, is to give globular sets $X(y,z)$ for all $y,z \in Z$ and maps $x_{y,z}:\{AX\}(y,z){\rightarrow}X(y,z)$. By lemma(\ref{lem:A-homs}) the $x_{y,z}$ amount to giving for each $k \in \N$ and $f:[k]_0{\rightarrow}X$ such that $f0=y$ and $fk=z$, a map
\[ x_f : \opA\limits_iX(f_{i-1},fi) \rightarrow X(y,z) \]
since $\opA\limits_iX(f_{i-1},fi)=\{Af^*X\}(0,k)$, that is $x_f=x_{y,z}\{A\overline{f}\}_{0,k}$. For $y,z \in Z$, one has a unique $f:[1]_0{\rightarrow}X$ given by $f0=y$ and $f1=z$. The naturality square for $\eta$ at $\overline{f}$ implies that $\{\eta_X\}_{y,z}=\{A\overline{f}\}_{0,1}\{\eta_{(X(y,z))}\}_{0,1}$ and the definition of $\overline{(\,\,)}$ says that $\{\eta_{(X(y,z))}\}_{0,1}=\overline{\eta}_{X(y,z)}$. Thus to say that a map $x:AX{\rightarrow}X$ satisfies the unit law of an $A$-algebra is to say that $x$ is the identity on $0$-cells and that the $x_f$ described above satisfy the unit axioms of an $\overline{A}$-category.

To say that $x$ satisfies the associative law is to say that for all $y,z \in Z$,
\begin{equation}\label{eq:assoc1} \xymatrix{{\{A^2X\}(y,z)} \ar[r]^-{\{\mu_X\}_{y,z}} \ar[d]_{\{Ax\}_{y,z}} & {\{AX\}(y,z)} \ar[d]^{x_{y,z}} \\ {\{AX\}(y,z)} \ar[r]_-{x_{y,z}} & {X(y,z)}} \end{equation}
commutes. Given $f:[m]_0{\rightarrow}X$ with $f0=y$ and $fm=z$, and $g:[k]_0{\rightarrow}Af^*X$ with $g0=0$ and $gk=m$, precomposing (\ref{eq:assoc1}) with the composite map
\begin{equation}\label{eq:copr} \xymatrix{{\{Ag^*Af^*X\}(0,k)} \ar[r]^-{\{A\overline{g}\}_{0,k}} & {\{A^2f^*X\}(0,m)} \ar[r]^-{\{A\overline{f}\}_{0,m}} & {\{A^2X\}(y,z)}} \end{equation}
and using lemma(\ref{lem:A-homs}) one can see that one obtains the commutativity of
\begin{equation}\label{eq:assoc2} \xymatrix{{\opA\limits_i\opA\limits_jX(f((i,j)-1),f(i,j))} \ar[r]^-{\overline{\mu}} \ar[d]_{\opA\limits_ix_{\{A\overline{f}\}g}} & {\opA\limits_{ij}X(f((i,j)-1),f(i,j))} \ar[d]^{x_f} \\ {\opA\limits_iX(g(i-1),gi)} \ar[r]_-{x_g} & {X(y,z)}} \end{equation}
where $1{\leq}i{\leq}k$, $1{\leq}j{\leq}m_i$, with the $m_i$ determined in the obvious way by $g$. That is, the associative law for $x$, namely (\ref{eq:assoc1}), implies the $\overline{A}$-category associative laws (\ref{eq:assoc2}). Conversely since the composites (\ref{eq:copr}) over all choices of $f$ and $g$ form a coproduct cocone by lemma(\ref{lem:A-homs}), (\ref{eq:assoc2}) also implies (\ref{eq:assoc1}). This completes the description of the object part of $\Alg A \iso \Enrich {\overline{A}}$.

Let $(X,x)$ and $(X',x')$ be $A$-algebras and $F_0:X_0{\rightarrow}X'_0$ be a function. To give $F:X{\rightarrow}X'$ with $0$-cell map $F_0$ is to give for all $y,z \in X_0$, maps $F_{y,z}:X(y,z){\rightarrow}X'(F_0y,F_0z)$. By lemma(\ref{lem:A-homs}) to say that $F$ is an algebra map is equivalent to saying that $F_0$ and the $F_{y,z}$ form an $\overline{A}$-functor. The isomorphism $\Alg A \iso \Enrich {\overline{A}}$ just described commutes with the forgetful functors into $\Set$ by definition.
\end{proof}
Early in the above proof we saw that $\overline{(\,\,)}$ sends morphisms in $\NColl {\ca T}$ to cartesian transformations. Since $\ca T^{\times}$ is tight by proposition(\ref{prop:wcat-tightpra}) and lemma(\ref{lem:Tcross-tight}), this implies by theorem(\ref{thm:bar-monoidal}) that $\overline{(\,\,)}$ may in fact be regarded as a strong monoidal functor
\[ \overline{(\,\,)} : \NColl {\ca T} \rightarrow \PraDist(\PSh {\G})/\ca T^{\times}. \]
For this manifestation of $\overline{(\,\,)}$ we have the following result.
\begin{theorem}\label{thm:bar-equiv}
The functor $\overline{(\,\,)}$ just described is an equivalence of categories $\NColl {\ca T} \catequiv \PraDist(\PSh {\G})/\ca T^{\times}$.
\end{theorem}
\begin{proof}
We will verify that $\overline{(\,\,)}$ is essentially surjective on objects and fully faithful. For a cartesian $\varepsilon:E{\rightarrow}\ca T^{\times}$ we now define $\alpha:A{\rightarrow}\ca T$ so that $\overline{\alpha}\iso\varepsilon$. For $X \in \PSh {\G}$ define $\{AX\}_0=X_0$, and for $x,y \in X_0$, define $\{AX\}(x,y)$ as a coproduct with coproduct injections
\[ c_f : \opE\limits_iX(f(i-1),fi) \rightarrow \{AX\}(x,y) \]
for each $f:[k]_0{\rightarrow}X$ with $f0=x$ and $fk=y$. This definition is functorial in $X$ in the obvious way. The components of $\alpha$ are identities on $0$-cells with the hom maps determined by the commutativity of
\begin{equation}\label{eq:def-alpha} \xymatrix{{\opE\limits_iX(f(i-1),fi)} \ar[r]^-{c_f} \ar[d]_{\varepsilon} & {\{AX\}(x,y)} \ar[d]^{\{\alpha_X\}_{x,y}} \\ {\{\ca Tf^*X\}(0,k)} \ar[r]_-{\{\ca T\overline{f}\}_{0,k}} & {\{\ca TX\}(x,y)}} \end{equation}
for all $f$ as above. Since $\PSh {\G}$ is extensive these squares are pullbacks, and so by lemma(\ref{lem:pb-hom}) $\alpha$ defined in this way is indeed cartesian. In the case where $X=(X_1,...,X_k)$ and $f$ is the identity on $0$-cells, one has $\{\ca T\overline{f}\}_{0,k}=\id$ and so (\ref{eq:def-alpha}) gives $\overline{\alpha} \iso \varepsilon$ as required. To verify fully faithfulness let $\alpha:A{\rightarrow}\ca T$ and $\beta:B{\rightarrow}\ca T$ be normalised collections, and $\phi:\overline{A}{\rightarrow}\overline{B}$ be a cartesian transformation. To finish the proof it suffices, by the tightness of $\ca T$ and $\ca T^{\times}$ and lemma(\ref{lem:tight->unique}), to define a cartesian transformation $\psi:A{\rightarrow}B$ unique such that $\overline{\psi}=\phi$. For $X \in \PSh {\G}$ and $f:[k]_0{\rightarrow}X$ this last equation says that such a $\psi$ must satisfy
\[ \{\psi_{f^*X}\}(0,k) = \phi_{X(f(i-1),fi)} : \{Af^*X\}(0,k) \rightarrow \{Bf^*X\}(0,k). \]
The cartesianness of $\psi$ and the tightness of $\ca T$ implies $\psi\beta=\alpha$ by lemma(\ref{lem:tight->unique}), and so $\{\psi_X\}_0$ is the identity. For $x,y \in X_0$ the map $\{\psi_X\}_{x,y}$ is determined by the commutativity of
\[ \xymatrix{{\{Af^*X\}(0,k)} \ar[r]^-{\{A\overline{f}\}_{0,k}} \ar[d]_{\{\psi_{f^*X}\}(0,k)} & {\{AX\}(x,y)} \ar[d]^{\{\psi_X\}_{x,y}} \\ {\{Bf^*X\}(0,k)} \ar[r]_-{\{B\overline{f}\}_{0,k}} & {\{BX\}(x,y)}} \]
for all $f$, since the $\{A\overline{f}\}_{0,k}$ form a coproduct cocone by lemma(\ref{lem:A-homs}). Note also that this square is a pullback by the extensivity of $\PSh {\G}$. This completes the definition of the components of $\psi$ and the proof that they are determined uniquely by $\phi$ and the equation $\overline{\psi}=\phi$, and so to finish the proof one must verify that the $\psi_X$ are cartesian natural in $X$. To this end let $F:X{\rightarrow}Y$. Since the components of $\alpha$ are identities in dimension $0$ it suffices by lemma(\ref{lem:pb-hom}) to show that for all $x,y \in X_0$ the squares
\begin{equation}\label{eq:natF} \xymatrix{{\{AX\}(x,y)} \ar[r]^-{\{AF\}_{x,y}} \ar[d]_{\{\psi_X\}_{x,y}} & {\{AY\}(F_0x,F_0y)} \ar[d]^{\{\psi_Y\}_{F_0x,F_0y}} \\ {\{BX\}(x,y)} \ar[r]_-{\{BF\}_{x,y}} & {\{BY\}(F_0x,F_0y)}} \end{equation}
are pullbacks. For all $f:[k]_0{\rightarrow}X$ one has $F\overline{f}=\overline{Ff}$ by definition, and so the composite square
\[ \xymatrix{{\{Af^*X\}(0,k)} \ar[r]^-{\{A\overline{f}\}_{0,k}} \ar[d]_{\{\psi_{f^*X}\}(0,k)} & {\{AX\}(x,y)} \ar[r]^-{\{AF\}_{x,y}} \ar[d]_{\{\psi_X\}_{x,y}} & {\{AY\}(F_0x,F_0y)} \ar[d]^{\{\psi_Y\}_{F_0x,F_0y}} \\ {\{Bf^*X\}(0,k)} \ar[r]_-{\{B\overline{f}\}_{0,k}} & {\{BX\}(x,y)} \ar[r]_-{\{BF\}_{x,y}} & {\{BY\}(F_0x,F_0y)}} \]
is a pullback, and so by the extensivity of $\PSh {\G}$ (\ref{eq:natF}) is indeed a pullback since the $\{A\overline{f}\}_{0,k}$ for all $f$ form a coproduct cocone.
\end{proof}
\begin{remark}\label{rem:eq-alt-view}
The equivalence of theorem(\ref{thm:bar-equiv}) could have been described differently. This alternative view involves the adjoint endofunctors $D$ and $\Sigma$ of $\PSh {\G}$. For $X \in \PSh {\G}$, $DX$ is obtained by discarding the $0$-cells and putting $\{DX\}_n=X_{n+1}$, $\Sigma{X}$ has one $0$-cell and $\Sigma{X}_{n+1}=X_n$ and one has $D \ladj \Sigma$. The effect of $D$ and $\Sigma$ on arrows provides an adjunction
\begin{equation}\label{eq:simple-equiv}
\xymatrix{{{\PSh {\G}}/\ca T1} \ar@<1.7ex>[rr]^-{D_{\ca T1}} \save \POS?="top" \restore && {\PSh {\G}/D\ca T1}
\ar@<1.7ex>[ll]^-{\Sigma_{D\ca T1}} \save \POS?="bot" \restore \POS "top"; "bot" **@{}; ?*{\perp}}, \end{equation}
and the right adjoint $\Sigma_{D\ca T1}$ is fully faithful since $\Sigma$ is. Thus (\ref{eq:simple-equiv}) restricts to an equivalence between the full subcategory $\ca N$ of $\PSh {\G}/\ca T1$ consisting of those $f:X{\rightarrow}\ca T1$ such that $X_0$ is singleton. Evaluating at $1$ gives an equivalence between $\NColl {\ca T}$ and the full subcategory of $\PSh {\G}/\ca T1$ just described. By evaluating at $1$ and by the definitions of $D$ and $\ca T1$ one obtains $\PraDist(\PSh {\G})/\ca T^{\times} \catequiv \PSh {\G}/D\ca T1$. Finally these equivalences fit together into a square
\[ \xymatrix{{\NColl {\ca T}} \ar[r]^-{\overline{(\,\,)}} \ar[d]_{\ev_1} & {\PraDist(\PSh {\G})/\ca T^{\times}} \ar[d]^{\ev_1} \\ {\ca N} \ar[r]_-{D_{\ca T1}} & {\PSh {\G}/D\ca T1}} \]
which one may easily verify commutes up to isomorphism. These equivalences $\ev_1$ really just express the equivalence of two different ways of viewing collections and their multitensorial analogues, and so modulo this, the equivalence from (\ref{eq:simple-equiv}) expresses in perhaps more concrete terms what $\overline{(\,\,)}$ does. However we have chosen to work with $\overline{(\,\,)}$ because this point of view makes clearer the relationship between algebras and enriched categories that we have expressed in theorem(\ref{thm:bar-monoidal}).
\end{remark}
Putting together theorem(\ref{thm:bar-equiv}) and theorem(\ref{thm:TMult-MTOp}) one obtains the equivalence between normalised $\ca T$-operads, $\ca T$-multitensors and $M\ca T$-operads.
\begin{corollary}\label{cor:3way-equiv}
$\NOp {\ca T} \catequiv \Mult {\ca T} \catequiv \Op {M\ca T}$.
\end{corollary}
%

%%%%%%%%%%%%%%%%%%%%%%%%%%%%%%%%%%%%%%%%%%%%%%%%%%
%%%%%%%%%%%%%%%%%%%%%%%%%%%%%%%%%%%%%%%%%%%%%%%%%%
\section{Finite dimensions and the algebras of $\ca T$}\label{sec:fin-dim}

We shall now explain how the results of this paper specialise to finite dimensions, and show how one can see that the algebras of $\ca T$ really are strict $\omega$-categories defined in the usual way by successive enrichment.

The category $\G_{{\leq}n}$ is defined to be the full subcategory of $\G$ consisting of the $k \in \N$ such that $0 \leq k \leq n$. The objects of $\PSh {\G}_{{\leq}n}$ are called $n$-globular sets. By definition the monad $\ca T$ on $\PSh {\G}$ restricts to $n$-globular sets: the description of $\ca TX_n$ depends only on the $k$-cells of $X$ for $k \leq n$. Thus one has a monad $\ca T_{{\leq}n}$ on $\PSh {\G}_{{\leq}n}$. Our description of $\ca T$ from section(\ref{sec:wcat-monad}) restricts also, and so the monads $\ca T_{{\leq}n}$ are p.r.a, coproduct preserving and tight. In fact, by direct inspection, everything we have done in this paper that has anything to do with $\ca T$ restricts to finite dimensions.

In particular for $n \in \N$, denoting by $\NColl {\ca T_{{\leq}1+n}}$ the category of normalised $(1+n)$-collections, whose objects are cartesian transformations $\alpha:A{\rightarrow}\ca T_{{\leq}1+n}$ whose components are identities in dimension $0$, one has a functor
\[ \overline{(\,\,)} : \NColl {\ca T_{{\leq}1+n}} \rightarrow \Dist(\PSh {\G}_{\leq n}) \]
whose object map is given by the formula
\[ \opA\limits_iX_i = \{AX\}(0,k) \]
where $A$ is a normalised $(1+n)$-collection, $k \in \N$ and $X_i \in \PSh {\G}_{{\leq}n}$ where $1 \leq i \leq k$, and $X \in \PSh {\G}_{{\leq}1+n}$ is defined as $X=(X_1,...,X_k)$. The finite dimensional analogue of theorem(\ref{thm:bar-monoidal}) is
\begin{theorem}\label{thm:bar-monoidal-finite}
The functor $\overline{(\,\,)}$ just described is a strong monoidal functor, and for a normalised $(1+n)$-operad $A$, one has an isomorphism $\Alg A \iso \Enrich {\overline{A}}$ commuting with the forgetful functors into $\Set$.
\end{theorem}
\noindent As before one may also regard $\overline{(\,\,)}$ as a strong monoidal functor
\[ \overline{(\,\,)} : \NColl {\ca T_{{\leq}1+n}} \rightarrow \PraDist(\PSh {\G}_{{\leq}n})/\ca T^{\times}_{{\leq}n}. \]
and the analogue of theorem(\ref{thm:bar-equiv}) is
\begin{theorem}\label{thm:bar-equiv-finite}
The functor $\overline{(\,\,)}$ just described is an equivalence of categories $\NColl {\ca T_{{\leq}1+n}} \catequiv \PraDist(\PSh {\G}_{{\leq}n})/\ca T^{\times}_{{\leq}n}$.
\end{theorem}
\noindent and so we have
\begin{corollary}\label{cor:3way-equiv-finite}
$\NOp {\ca T_{{\leq}1+n}} \catequiv \Mult {\ca T_{{\leq}n}} \catequiv \Op {M\ca T_{{\leq}n}}$.
\end{corollary}
\noindent One can think of $n$ as an ordinal instead of a natural number, and then the original results from section(\ref{ssec:corr}) correspond to the case $n=\omega$.

All along we have been working with the monads $\ca T_{{\leq}n}$ as formally defined combinatorial objects. Given the results of this paper however, it is now easy to see that their algebras are indeed strict $n$-categories. The usual definition of strict $n$-categories is by successive enrichment. One defines $\Enrich 0 = \Set$ and $\Enrich {(1+n)}=\Enrich {(\Enrich n)}$ for $n \in \N$ where $\Enrich n$ is regarded as monoidal via cartesian product. Recasting this a little more formally, $\Enrich {(-)}$ is an endofunctor of the full subcategory of $\CAT$ consisting of categories with finite products. Writing $0$ for the terminal object of this category, that is the terminal category, one has by functoriality a sequence
\[ \xygraph{!{0;(1.5,0):}
0:@{<-}[l]{\Enrich 0}:@{<-}[l]{\Enrich 1}:@{<-}[l]{\Enrich 2}:@{<-}[l]{\Enrich 3}:@{<-}[l(.7)]{...}} \]
Explicitly the maps in this diagram are the obvious forgetful functors. The limit of this diagram is formed as in $\CAT$, and provides the definition of the category $\Enrich {\omega}$. Then by theorem(\ref{thm:bar-monoidal}) and proposition(\ref{prop:Tcross-cats}) we have isomorphisms
\[ \phi_{n} : \Alg {\ca T_{{\leq}1+n}} \rightarrow \Enrich {(\Alg {\ca T_{{\leq}n}})} \]
Let us write $\Enr$ for the endofunctor $\ca V \mapsto \Enrich {\ca V}$ that we have just been considering. The isomorphisms $\phi_n$ are natural in the sense of the following lemma, which enables us to then formally identify the algebras of $\ca T$ in theorem(\ref{thm:wcat-Talg}).
\begin{lemma}\label{lem:nat-phi}
For $n \in \N$ let $\tr_n:\Alg {\ca T_{{\leq}1+n}}{\rightarrow}\Alg {\ca T_{{\leq}n}}$ be the forgetful functor given by truncation. The square
\[ \xymatrix{{\Alg {\ca T_{{\leq}2+n}}} \ar[r]^-{\tr_{1+n}} \ar[d]_{\phi_{1+n}} & {\Alg {\ca T_{{\leq}1+n}}} \ar[d]^{\phi_n} \\ {\Enrich {(\Alg {\ca T_{{\leq}1+n}})}} \ar[r]_-{\Enr(\tr_n)} & {\Enrich {(\Alg {\ca T_{{\leq}n}})}}} \]
commutes for all $n \in \N$.
\end{lemma}
\begin{proof}
One obtains $\phi_n$ explicitly as the composite of two isomorphisms
\[ \Alg {\ca T_{{\leq}1+n}} \rightarrow \Enrich {\ca T^{\times}_{{\leq}n}} \rightarrow \Enrich {(\Alg {\ca T_{{\leq}n}})} \]
the first of which is described explicitly in the proof of theorem(\ref{thm:bar-monoidal}), and the second in the proof of proposition(\ref{prop:Tcross-cats}), and using these descriptions one may easily verify directly the desired naturality.
\end{proof}
\begin{theorem}\label{thm:wcat-Talg}
For $0 \leq n \leq \omega$, $\Alg {\ca T_{{\leq}n}} \iso \Enrich n$.
\end{theorem}
\begin{proof}
Write $t:\Enrich 0{\rightarrow}0$ for the unique functor. By the definition of $\Enrich {\omega}$ it suffices to provide isomorphisms $\psi_n:\Alg {\ca T_{{\leq}n}}{\rightarrow}\Enrich n$ for $n \in \N$ natural in the sense that 
\[ \xymatrix{{\Alg {\ca T_{{\leq}1+n}}} \ar[r]^-{\tr_n} \ar[d]_{\psi_{1+n}} & {\Alg {\ca T_{{\leq}n}}} \ar[d]^{\psi_n} \\ {\Enrich {(1+n)}} \ar[r]_-{\Enr^{1+n}(t)} & {\Enrich n}} \]
commutes for all $n$. Take $\psi_0=1_{\Set}$ and by induction define $\psi_{1+n}$ as the composite
\[ \xymatrix{{\Alg {\ca T_{{\leq}1+n}}} \ar[r]^-{\phi_n} & {\Enrich {(\Alg {\ca T_{{\leq}n}})}} \ar[r]^-{\Enr(\psi_n)} & {\Enrich {(1+n)}}}. \]
The case $n=0$ for $\psi$'s naturality comes from the fact that the isomorphisms that comprise $\phi_1$ (see lemma(\ref{lem:nat-phi})) are defined over $\Set$. The inductive step follows easily from lemma(\ref{lem:nat-phi}).
\end{proof}
%

%%%%%%%%%%%%%%%%%%%%%%%%%%%%%%%%%%%%%%%%%%%%%%%%%%
%%%%%%%%%%%%%%%%%%%%%%%%%%%%%%%%%%%%%%%%%%%%%%%%%%
\section{Acknowledgements}

The first author gratefully acknowledges for the financial support of Scott Russell Johnson Memorial Foundation. Both authors are also grateful for the financial support of  Max Planck Institut f\"{u}r Mathematik and the Australian Research Council grant No.~DP0558372. This paper was completed while the second author was a postdoc at Macquarie University in Sydney Australia and at the PPS lab in Paris, and he would like to thank these institutions for their hospitality and pleasant working conditions. We would both like to acknowledge the hospitality of the Max Planck Institute where some of this work was carried out. Moreover we are also indebted to Clemens Berger, Denis-Charles Cisinski and Paul-Andr\'e Melli\`es for interesting discussions on the substance of this paper.

\bibliographystyle{plain}
\bibliography{mtens}

\end{document}